\documentclass[leqno,draft]{article}

\newtheorem{theorem}{Theorem}
\newtheorem{lemma}[theorem]{Lemma}
\newtheorem{proposition}[theorem]{Proposition}
\newtheorem{definition}[theorem]{Definition}
\newtheorem{corollary}[theorem]{Corollary}

\newcommand{\begintheorem}{\addtocounter{equation}{1}\begin{theorem}}
\newcommand{\beginlemma}{\addtocounter{equation}{1}\begin{lemma}}
\newcommand{\beginproposition}{\addtocounter{equation}{1}\begin{proposition}}
\newcommand{\begindefinition}{\addtocounter{equation}{1}\begin{definition}}
\newcommand{\begincorollary}{\addtocounter{equation}{1}\begin{corollary}}

\begin{document}

\title{$p$-Adic Heisenberg Cantor sets}

\author{Stephen Semmes \\
        Rice University}

\date{}

\maketitle

\begin{abstract}
These informal notes deal with $p$-adic versions of Heisenberg groups
and related matters.
\end{abstract}

\tableofcontents

\section{Ultrametric spaces}
\label{ultrametric spaces}
\setcounter{equation}{0}

        Let $(M, d(x, y))$ be a metric space.  Thus $d(x, y)$ is a
nonnegative real-valued function defined for $x, y \in M$ such that
$d(x, y) = 0$ if and only if $x = y$,
\begin{equation}
\label{d(x, y) = d(y, x)}
        d(x, y) = d(y, x)
\end{equation}
for every $x, y \in M$, and
\begin{equation}
\label{d(x, z) le d(x, y) + d(y, z)}
        d(x, z) \le d(x, y) + d(y, z)
\end{equation}
for every $x, y, z \in M$.  As usual, this last condition is known as
the triangle inequality.  If
\begin{equation}
\label{d(x, z) le max (d(x, y), d(y, z))}
        d(x, z) \le \max (d(x, y), d(y, z))
\end{equation}
for every $x, y, z \in M$, then $d(x, y)$ is said to be an
\emph{ultrametric} on $M$.  Of course, (\ref{d(x, z) le max (d(x, y),
d(y, z))}) implies the ordinary triangle inequality (\ref{d(x, z) le
d(x, y) + d(y, z)}).

        If $(M, d(x, y))$ is a metric space, $x \in M$, and $r > 0$,
then the open ball in $M$ centered at $x$ with radius $r$ is defined
as usual by
\begin{equation}
\label{B(x, r) = {y in M : d(x, y) < r}}
        B(x, r) = \{y \in M : d(x, y) < r\},
\end{equation}
and the closed ball centered at $x$ with radius $r \ge 0$ is defined by
\begin{equation}
\label{overline{B}(x, r) = {y in M : d(x, y) le r}}
        \overline{B}(x, r) = \{y \in M : d(x, y) \le r\}.
\end{equation}
In any metric space, open balls are open sets, and closed balls are
closed sets.  In an ultrametric space, open balls are also closed
sets, and closed balls are also open sets.  This is because
\begin{equation}
\label{B(w, r) subseteq M backslash B(x, r)}
        B(w, r) \subseteq M \backslash B(x, r)
\end{equation}
for every $w \in M \backslash B(x, r)$, and
\begin{equation}
\label{overline{B}(z, r) subseteq overline{B}(x, r)}
        \overline{B}(z, r) \subseteq \overline{B}(x, r)
\end{equation}
for every $z \in \overline{B}(x, r)$, by the ultrametric version of
the triangle inequality.  In particular, ultrametric spaces are
totally disconnected.

        As a class of examples, let $A$ be a finite set with at least
two elements, let $\rho$ be a positive real number with $\rho < 1$,
and let $X$ be the set of sequences $x = \{x_j\}_{j = 1}^\infty$ such
that $x_j \in A$ for each $j$.  If $x, y \in X$ and $x \ne y$, then
let $n(x, y)$ be the largest nonnegative integer such that
\begin{equation}
\label{x_j = y_j for j le n(x, y)}
        x_j = y_j \quad\hbox{for } j \le n(x, y),
\end{equation}
and put
\begin{equation}
\label{d_rho(x, y) = rho^{-n(x, y)}}
        d_\rho(x, y) = \rho^{-n(x, y)}.
\end{equation}
If we put $d_\rho(x, x) = 0$ for each $x \in X$, then one can check
that this defines an ultrametric on $X$.  Equivalently, $X$ is the
same as the Cartesian product of a sequence of copies of $A$, which
leads to a product topology on $X$, associated to the discrete
topology on each copy of $A$.  Because $A$ is finite, it is compact,
and hence $X$ is compact with respect to the product topology.  It is
easy to see that the topology on $X$ determined by $d_\rho(x, y)$ is
the same as the product topology on $X$ for each $\rho \in (0, 1)$.
If instead we took $\rho = 1$ in this construction, then $d_\rho(x,
y)$ would be the discrete metric on $X$, which corresponds to the
discrete topology on $X$.

        Let $(M, d(x, y))$ be a metric space, and let $\{x_n\}_{n =
1}^\infty$ be a sequence of elements of $M$.  If $\{x_n\}_{n =
1}^\infty$ is a Cauchy sequence, then
\begin{equation}
        \lim_{n \to \infty} d(x_n, x_{n + 1}) = 0.
\end{equation}
Of course, the converse does not work in general, even for the real
line with the standard metric.  On the real line, the converse would
imply that an infinite series $\sum_{j = 1}^\infty a_j$ of real
numbers converges when $\lim_{j \to \infty} a_j = 0$.  However, it is
easy to see that the converse does hold when $d(x, y)$ is an
ultrametric on $M$.

\section{$p$-Adic numbers}
\label{p-adic numbers}
\setcounter{equation}{0}

        Let $p = 2, 3, 5, 7, \ldots$ be a prime number.  The
\emph{$p$-adic absolute value} $|x|_p$ of a rational number $x$ is
defined as follows.  If $x = 0$, then $|x|_p = 0$.  Otherwise, if $x
\ne 0$, then $x$ can be expressed as
\begin{equation}
        x = p^j \, \frac{a}{b},
\end{equation}
where $a$, $b$, and $j$ are integers, and $a$, $b$ are nonzero and not
divisible by $p$.  In this case, we put
\begin{equation}
\label{|x|_p = p^{-j}}
        |x|_p = p^{-j}.
\end{equation}
In particular, $|x|_p$ is equal to $0$ or to an integer power of $p$ for
every $x$ in the set ${\bf Q}$ of rational numbers.  It is easy to check that
\begin{equation}
\label{|x + y|_p le max(|x|_p, |y|_p)}
        |x + y|_p \le \max (|x|_p, |y|_p)
\end{equation}
and
\begin{equation}
\label{|x y|_p = |x|_p |y|_p}
        |x \, y|_p = |x|_p \, |y|_p
\end{equation}
for every $x, y \in {\bf Q}_p$.

        The \emph{$p$-adic metric} is defined on ${\bf Q}$ by
\begin{equation}
\label{d_p(x, y) = |x - y|_p}
        d_p(x, y) = |x - y|_p.
\end{equation}
This defines an ultrametric on ${\bf Q}$, because of (\ref{|x + y|_p
le max(|x|_p, |y|_p)}).  However, ${\bf Q}$ is not complete with
respect to this ultrametric.  The field ${\bf Q}_p$ of $p$-adic
numbers is obtained by completing ${\bf Q}$ with respect to $d_p(x,
y)$, in the same way that the real numbers can be obtained by
completing the rational numbers with respect to the standard metric.
More precisely, the usual arithmetic operations of addition and
multiplication are extended from ${\bf Q}$ to ${\bf Q}_p$, as well as
the $p$-adic absolute value and metric.  The $p$-adic absolute value
of any element of ${\bf Q}_p$ is either $0$ or an integer power of
$p$, and the extension of the $p$-adic absolute value to ${\bf Q}_p$
still satisfies (\ref{|x + y|_p le max(|x|_p, |y|_p)}) and (\ref{|x
y|_p = |x|_p |y|_p}).  The $p$-adic metric on ${\bf Q}_p$ is still
expressed in terms of the $p$-adic absolute value on ${\bf Q}_p$ as in
(\ref{d_p(x, y) = |x - y|_p}), and ${\bf Q}$ is dense in ${\bf Q}_p$
with respect to the $p$-adic metric.

        As in the context of real numbers, an infinite series
\begin{equation}
\label{sum_{j = 1}^infty a_j}
        \sum_{j = 1}^\infty a_j
\end{equation}
of $p$-adic numbers is said to converge in ${\bf Q}_p$ if and only if
the corresponding sequence of partial sums
\begin{equation}
\label{s_n = sum_{j = 1}^n a_j}
        s_n = \sum_{j = 1}^n a_j
\end{equation}
converges in ${\bf Q}_p$.  Because of completeness, this happens if
and only if $\{s_n\}_{n = 1}^\infty$ is a Cauchy sequence in ${\bf
Q}_p$.  If $\{s_n\}_{n = 1}^\infty$ is a Cauchy sequence in ${\bf
Q}_p$, then
\begin{equation}
\label{lim_{j to infty} a_j = 0}
        \lim_{j \to \infty} a_j = 0,
\end{equation}
as in the case of real numbers.  However, in contrast with the
situation for real numbers, (\ref{lim_{j to infty} a_j = 0}) implies
that $\{s_n\}_{n = 1}^\infty$ is a Cauchy sequence in ${\bf Q}_p$, and
hence converges, because the $p$-adic metric is an ultrametric.

\section{$p$-Adic integers}
\label{$p$-adic integers}
\setcounter{equation}{0}

        Let ${\bf Z}$ be the set of integers, which can also be
considered as a subring of ${\bf Q}$.  If $x \in {\bf Z}$, then $|x|_p
\le 1$, by definition of the $p$-adic absolute value.  Put
\begin{equation}
        {\bf Z}_p = \{x \in {\bf Q}_p : |x|_p \le 1\},
\end{equation}
so that ${\bf Z} \subseteq {\bf Z}_p$.  The elements of ${\bf Z}_p$
are called $p$-adic integers, and it is easy to see that they form a
subring of ${\bf Q}_p$, by (\ref{|x + y|_p le max(|x|_p, |y|_p)}) and
(\ref{|x y|_p = |x|_p |y|_p}).

        Let $x$ be a nonzero rational number, which can be expressed
as $p^j \, (a/b)$, where $a$, $b$, $j$ are integers, $a, b \ne 0$, and
$a$, $b$ are not divisible by $p$, as in the previous section.
Because the ring ${\bf Z} / p \, {\bf Z}$ of integers modulo $p$ form
a field, there is a nonzero integer $c$ such that $b \, c \equiv 1$
modulo $p$.  This implies that
\begin{equation}
        x = p^j \, \frac{a \, c}{b \, c} = p^j \, \frac{a \, c}{1 - p \, e}
\end{equation}
for some $e \in {\bf Z}$.  Note that
\begin{equation}
\label{sum_{l = 0}^n (p e)^l = frac{1 - (p e)^{n + 1}}{1 - p e}}
        \sum_{l = 0}^n (p \, e)^l = \frac{1 - (p \, e)^{n + 1}}{1 - p \, e}
\end{equation}
for each positive integer $n$, when $(p \, e)^l$ is interpreted as being
equal to $1$ when $l = 0$, even when $e = 0$.  Hence $\sum_{l = 0}^\infty
(p \, e)^l$ converges with respect to the $p$-adic metric, since
\begin{equation}
        |p \, e|_p = |p|_p \, |e|_p \le \frac{1}{p} < 1,
\end{equation}
and in fact
\begin{equation}
        \sum_{l = 0}^\infty (p \, e)^l = \frac{1}{1 - p \, e}.
\end{equation}
In particular, it follows that $1/(1 - p \, e)$ can be approximated by
integers with respect to the $p$-adic metric.  If $|x|_p \le 1$, so
that $j \ge 0$, then we also get that $x$ can be approximated by
integers with respect to the $p$-adic metric.

        Suppose now that $x \in {\bf Z}_p$.  Of course, $x$ can be
approximated by rational numbers with respect to the $p$-adic metric,
since ${\bf Q}$ is dense in ${\bf Q}_p$.  If $y \in {\bf Q}$ and
$|x - y|_p \le 1$, then
\begin{equation}
        |y|_p \le \max (|x|_p, |x - y|_p) \le 1.
\end{equation}
This implies that $y$ can be approximated by ordinary integers with
respect to the $p$-adic metric, by the previous argument.  Thus
$p$-adic integers can be approximated by ordinary integers with
respect to the $p$-adic metric.  Note that ${\bf Z}_p$ is a closed set
with respect to the $p$-adic metric, because it is the same as the
closed ball in ${\bf Q}_p$ with center $0$ and radius $1$.  It follows that
\begin{equation}
        {\bf Z}_p = \overline{\bf Z} \quad\hbox{in } {\bf Q}_p,
\end{equation}
which is to say that ${\bf Z}_p$ is the same as the closure of ${\bf
Z}$ in ${\bf Q}_p$ with respect to the $p$-adic metric, since ${\bf Z}
\subseteq {\bf Z}_p$.

\section{Ideals and quotients}
\label{ideals, quotients}
\setcounter{equation}{0}

        Let $j$ be a nonnegatuve integer, and consider
\begin{equation}
        p^j \, {\bf Z}_p = \{p^j \, x : x \in {\bf Z}_p\}
                         = \{y \in {\bf Q}_p : |y|_p \le p^{-j}\}.
\end{equation}
This is a closed subset of ${\bf Q}_p$, and also an ideal in ${\bf
Z}_p$.  Of course,
\begin{equation}
        p^j \, {\bf Z} = \{p^j \, x : x \in {\bf Z}\}
\end{equation}
is also an ideal in ${\bf Z}$, which is contained in $p^j \, {\bf Z}_p$.
It follows from the discussion in the previous section that $p^j \, {\bf Z}_p$
is the same as the closure of $p^j \, {\bf Z}$ in ${\bf Q}_p$.  Moreover,
\begin{equation}
\label{p^j {bf Z} = {bf Z} cap (p^j {bf Z}_p)}
        p^j \, {\bf Z} = {\bf Z} \cap (p^j \, {\bf Z}_p).
\end{equation}

        Because $p^j \, {\bf Z}_p$ is an ideal in ${\bf Z}_p$, the
quotient ${\bf Z}_p / p^j {\bf Z}_p$ is a commutative ring in a
natural way.  The obvious inclusion of ${\bf Z}$ in ${\bf Z}_p$
leads to a mapping
\begin{equation}
\label{{bf Z} mapsto {bf Z}_p p^j {bf Z}_p}
        {\bf Z} \mapsto {\bf Z}_p / p^j \, {\bf Z}_p,
\end{equation}
which is a homomorphism between commutative rings.  The kernel of this
homomorphism is equal to $p_j \, {\bf Z}$, by (\ref{p^j {bf Z} = {bf
Z} cap (p^j {bf Z}_p)}), and so we get an injective homomorphism
\begin{equation}
\label{{bf Z} / p^j {bf Z} mapsto {bf Z}_p / p^j {bf Z}_p}
        {\bf Z} / p^j \, {\bf Z} \mapsto {\bf Z}_p / p^j \, {\bf Z}_p.
\end{equation}

        If $x \in {\bf Z}_p$, then there is a $y \in {\bf Z}$ such that
\begin{equation}
        |x - y|_p \le p^{-j},
\end{equation}
because ${\bf Z}$ is dense in ${\bf Z}_p$.  Equivalently, $x - y \in
p^j \, {\bf Z}_p$.  This implies that the mapping (\ref{{bf Z} / p^j
{bf Z} mapsto {bf Z}_p / p^j {bf Z}_p}) is a surjection, and hence an
isomorphism.  Thus
\begin{equation}
        {\bf Z}_p / p^j \, {\bf Z}_p \cong {\bf Z} / p^j \, {\bf Z}.
\end{equation}

        In particular, ${\bf Z}_p / p^j \, {\bf Z}_p$ is a finite set
with $p^j$ elements for each nonnegative integer $j$.  This shows that
${\bf Z}_p$, which is the closed unit ball in ${\bf Q}_p$, can be
expressed as the union of $p^j$ translates of $p^j \, {\bf Z}_p$,
which is a union of $p^j$ closed balls of radius $p^{-j}$ in ${\bf
Q}_p$.  It follows that ${\bf Z}_p$ is totally bounded in ${\bf Q}_p$,
and is therefore compact, since ${\bf Z}_p$ is also closed, and ${\bf
Q}_p$ is complete.

\section{Heisenberg groups}
\label{heisenberg groups}
\setcounter{equation}{0}

        Let $R$ be a commutative ring, not necessarily with a
multiplicative identity element.  If $n$ is a positive integer, then
the set $R^n$ of $n$-tuples $x = (x_1, \ldots, x_n)$ of elements of
$R$ is also a commutative ring with respect to coordinatewise addition
and multiplication.  However, we shall be more concerned with $R^n$ as
a module over $R$, with respect to coordinatewise addition and
``scalar multiplication'', defined by
\begin{equation}
        r \, x = (r \, x_1, \ldots, r \, x_n)
\end{equation}
for each $r \in R$ and $x \in R^n$.  In particular, if $R$ is a field,
then $R^n$ is a vector space over $R$.

        Similarly, $R^n \times R^n$ is a module over $R$, and a vector
space when $R$ is a field.  Let us use the notation $z = (x, y)$ with
$x, y \in R^n$ for elements of $R^n \times R^n$, and $r \, z = (r \,
x, r \, y)$ for scalar multiplication with $r \in R$.  Consider
\begin{equation}
        B(z, z') = \sum_{j = 1}^n (x_j \, y'_j - x'_j \, y_j)
\end{equation}
for $z = (x, y), z' = (x', y') \in R^n \times R^n$.  Note that this is
linear in $z$ and $z'$ separately, with respect to coordinatewise
addition and scalar multiplication by elements of $R$.  Of course,
$B(z, z')$ is also antisymmetric in $z$ and $z'$, in the sense that
\begin{equation}
        B(z', z) = - B(z, z').
\end{equation}

        The $n$th Heisenberg group over $R$ is initially
defined as a set by
\begin{equation}
        H_n(R) = R^n \times R^n \times R \cong (R^n \times R^n) \times R.
\end{equation}
Using the latter identification, elements of $H_n(R)$ will often be
denoted $(z, t)$, where $z = (x, y) \in R^n \times R^n$ and $t \in R$.
The group operation on $H_n(R)$ is defined by
\begin{equation}
        (z, t) \diamond (z', t') = (z + z', t + t' + B(z, z')).
\end{equation}
It is easy to see that this is associative, because
\begin{eqnarray}
\lefteqn{((z, t) \diamond (z', t')) \diamond (z'', t'')
         =  (z + z', t + t' + B(z, z')) \diamond (z'', t'')} \\
 & = & (z + z' + z'', t + t' + t'' + B(z, z') + B(z + z', z'')) \nonumber
\end{eqnarray}
is equal to
\begin{eqnarray}
\lefteqn{(z, t) \diamond ((z', t') \diamond (z'', t''))
         = (z, t) \diamond (z' + z'', t' + t'' + B(z', z''))} \\
 & = & (z + z' + z'', t + t' + t'' + B(z, z' + z'') + B(z', z'')), \nonumber
\end{eqnarray}
since
\begin{eqnarray}
        B(z, z') + B(z + z', z'') & = & B(z, z') + B(z, z'') + B(z', z'') \\
                            & = & B(z, z' + z'') + B(z', z''). \nonumber
\end{eqnarray}
It is easy to see that
\begin{equation}
        (0, 0) \diamond (z, t) = (z, t) \diamond (0, 0) = (z, t)
\end{equation}
and
\begin{equation}
        (-z, -t) \diamond (z, t) = (z, t) \diamond (-z, -t) = (0, 0)
\end{equation}
for every $(z, t) \in H_n(R)$.  Thus $(0, 0)$ is the identity element
in $H_n(R)$, and $(-z, -t)$ is the inverse of $(z, t)$ in $H_n(R)$.
This shows that $H_n(R)$ is a group.

        Of course, $H_n(R)$ is not normally commutative, because $B(z,
z')$ is anti-symmetric in $z$, $z'$.  However, it could be that $R$
has ``characteristic $2$'', in the sense that
\begin{equation}
        r + r = 0
\end{equation}
for every $r \in R$, so that the additive inverse of $r$ is equal to
$r$ for every $r \in R$, which is to say that $r = -r$ for every $r \in R$.
In this case, $B(z, z')$ would also be symmetric in $z$, $z'$,
and $H_n(R)$ would be commutative.

\section{Some variants}
\label{some variants}
\setcounter{equation}{0}

        Let $A$ be an abelian group, in which the group operation is
expressed by addition.  If $N$ is a positive integer, then the set
$A^N$ of $N$-tuples of elements of $A$ is also an abelian group, with
respect to coordinatewise addition.  Suppose that $B(w, z)$ is an
$A$-valued function defined for $w, z \in A^N$ which is additive in
each variable, in the sense that
\begin{equation}
        B(w + w', z) = B(w, z) + B(w', z)
\end{equation}
and
\begin{equation}
        B(w, z + z') = B(w, z) + B(w, z')
\end{equation}
for every $w, w', z, z' \in A^N$.  This includes the situation
described in the previous section, with $A = R$ as a group with
respect to addition, and $N = 2 n$.

        Consider the binary operation defined on $A^N \times A$ by
\begin{equation}
\label{(w, s) diamond (z, t) = (w + z, s + t + B(w, z))}
        (w, s) \diamond (z, t) = (w + z, s + t + B(w, z)).
\end{equation}
Again, this includes the group operation in the previous section as a
special case.  It is easy to see that (\ref{(w, s) diamond (z, t) = (w
+ z, s + t + B(w, z))}) is associative, for exactly the same reasons
as before.  If $0$ is the additive identity element in $A$, then $(0,
\ldots, 0)$ is the additive identity element $(0, \ldots, 0)$ in
$A^N$, which is also denoted $0$.  If $w$, $z$ are elements of $A^N$, then
\begin{equation}
        B(w, z) = B(w + 0, z) = B(w, z) + B(0, z),
\end{equation}
and hence $B(0, z) = 0$.  Similarly, $B(w, 0) = 0$, and it follows that
\begin{equation}
        (0, 0) \diamond (z, t) = (z, t) \diamond (0, 0) = (z, t)
\end{equation}
for every $(z, t) \in A^N \times A$.  Thus $(0, 0)$ is the identity
element for (\ref{(w, s) diamond (z, t) = (w + z, s + t + B(w, z))}),
and we would like to check that every element of $A^N \times A$ has an
inverse with respect to (\ref{(w, s) diamond (z, t) = (w + z, s + t +
B(w, z))}).

        Observe that
\begin{equation}
        B(w, z) + B(-w, z) = B(0, z) = 0
\end{equation}
and
\begin{equation}
        B(w, z) + B(w, -z) = B(w, 0) = 0
\end{equation}
for every $w, z \in A^N$.  If $B(z, z) = 0$ for every $z \in A^N$, so
that $B(z, -z) = 0$ for every $z \in A^N$ too, then
\begin{equation}
\label{(z, t) diamond (-z, -t) = (-z, -t) diamond (z, t) = (0, 0)}
        (z, t) \diamond (-z, -t) = (-z, -t) \diamond (z, t) = (0, 0)
\end{equation}
for every $(z, t) \in A^N \times A$.  In this case, $(-z, -t)$ is the
inverse of $(z, t)$ with respect to (\ref{(w, s) diamond (z, t) = (w +
z, s + t + B(w, z))}).  Otherwise,
\begin{equation}
\label{(z, t) diamond (-z, -t + B(z, z)) = ... = (0, 0)}
        (z, t) \diamond (-z, -t + B(z, z))
           = (-z, -t + B(z, z)) \diamond (z, t) = (0, 0)
\end{equation}
for every $(z, t) \in A^N \times A$, and so $(-z, -t + B(z, z))$ is
the inverse of $(z, t)$ with respect to (\ref{(w, s) diamond (z, t) =
(w + z, s + t + B(w, z))}).

        Note that
\begin{equation}
        B(w + z, w + z) = B(w, w) + B(w, z) + B(z, w) + B(z, z)
\end{equation}
for every $w, z \in A^N$.  If $B(z, z) = 0$ for each $z \in A^N$, then
it follows that
\begin{equation}
\label{B(w, z) + B(z, w) = 0}
        B(w, z) + B(z, w) = 0
\end{equation}
for every $w, z \in A^N$.  If $A$ is a field with characteristic not
equal to $2$, for example, then (\ref{B(w, z) + B(z, w) = 0}) implies
that $B(z, z) = 0$ for every $z \in A^N$.  However, if $A$ is a field
with characteristic $2$, for instance, then (\ref{B(w, z) + B(z, w) =
0}) is the same as saying that $B(w, z)$ is symmetric in $w$ and $z$,
and this does not necessarily imply that $B(z, z) = 0$ for each $z \in
A^N$.

\section{Changing variables}
\label{changing variables}
\setcounter{equation}{0}

        Let us continue with the same notations and hypotheses as in
the previous section.  Let $C(w, z)$ be another $A$-valued function of
$w, z \in A^N$ which is additive in each variable, and put
\begin{equation}
        \widetilde{B}(w, z) = B(w, z) + C(w, z) + C(z, w).
\end{equation}
This leads to another group operation on $A^N \times A$, defined by
\begin{equation}
 (w, s) \widetilde{\diamond} (z, t) = (w + z, s + t + \widetilde{B}(w, z)),
\end{equation}
as in (\ref{(w, s) diamond (z, t) = (w + z, s + t + B(w, z))}).

        Consider the mapping $\phi : A^N \times A \to A^N \times A$, defined by
\begin{equation}
\label{phi((z, t)) = (z, t + C(z, z))}
        \phi((z, t)) = (z, t + C(z, z)).
\end{equation}
It is easy to see that $\phi$ is a one-to-one mapping from $A^N \times
A$ onto itself, and that $\phi((0, 0)) = (0, 0)$.  Moreover,
\begin{eqnarray}
\label{phi((w, s) diamond (z, t)) = ... = (w + z, s + t + B(w, z) + C(w + z))}
        \phi((w, s) \diamond (z, t)) & = & \phi((w + z, s + t + B(w, z))) \\
            & = & (w + z, s + t + B(w, z) + C(w + z, w + z)) \nonumber
\end{eqnarray}
for every $(w, s), (z, t) \in A^N \times A$.  Of course,
\begin{equation}
        C(w + z, w + z) = C(w, w) + C(w, z) + C(z, w) + C(z, z),
\end{equation}
as in the previous section, and so we get that
\begin{equation}
 \phi((w, s) \diamond (z, t)) = \phi((w, s)) \widetilde{\diamond} \phi((z, t))
\end{equation}
for every $(w, s), (z, t) \in A^N \times A$.  Thus $\phi$ defines an
isomorphism from $A^N \times A$ equipped with $\diamond$ as the group
operation onto $A^N \times A$ equipped with $\widetilde{\diamond}$ as
the group operation.

        If $B(w, z) = 0$ for every $w, z \in A^N$, then
\begin{equation}
 (w, s) \widetilde{\diamond} (z, t) = (w + z, s + t + C(w, z) + C(z, w))
                                    = (z, t) \widetilde{\diamond} (w, s).
\end{equation}
In this case, $\phi$ defines an isomorphism from $A^N \times A$ as a
commutative group with respect to coordinatewise addition onto $A^N
\times A$ equipped with $\widetilde{\diamond}$.

\section{Subgroups and quotients}
\label{subgroups, quotients}
\setcounter{equation}{0}

        Let us continue with the same notations and hypotheses as in
the previous sections.  Suppose that $A_0$ is a subgroup of $A$, so
that $A_0^N \times A_0$ is a subgroup of $A^N \times A$ as an abelian
group with respect to coordinatewise addition.  If we also have that
\begin{equation}
        B(w, z) \in A_0
\end{equation}
for every $w, z \in A_0$, then we get that
\begin{equation}
 (w, s) \diamond (z, t) = (w + z, s + t + B(w, z)) \in A_0^N \times A_0
\end{equation}
for every $(w, s), (z, t) \in A_0^N \times A_0$.  Similarly, the inverse
\begin{equation}
        (w, s)^{-1} = (-w, -s + B(w, w))
\end{equation}
of $(w, s)$ with respect to (\ref{(w, s) diamond (z, t) = (w + z, s +
t + B(w, z))}) is an element of $A_0^N \times A_0$ for every $(w, s)$
in $A_0^N \times A_0$.  This shows that $A_0^N \times A_0$ is a
subgroup of $A^N \times A$ with respect to (\ref{(w, s) diamond (z, t)
= (w + z, s + t + B(w, z))}) under these conditions.

        Observe that
\begin{eqnarray}
\label{((w, s) diamond (z, t)) diamond (w, s)^{-1} = ...}
\lefteqn{((w, s) \diamond (z, t)) \diamond (w, s)^{-1}} \\
         & = & (w + z, s + t + B(w, z)) \diamond (-w, -s + B(w, w)) \nonumber\\
         & = & (z, t + B(w, z) + B(w, w) + B(w + z, -w)) \nonumber\\
         & = & (z, t + B(w, z) - B(z, w)) \nonumber
\end{eqnarray}
for every $(w, s), (z, t) \in A^N \times A$.  If
\begin{equation}
        B(w, z) - B(z, w) \in A_0
\end{equation}
for every $w \in A^N$ and $z \in A_0^N$, then it follows that $A_0^N
\times A_0$ is a normal subgroup of $A^N \times A$ with respect to
(\ref{(w, s) diamond (z, t) = (w + z, s + t + B(w, z))}).  In
particular, this holds when
\begin{equation}
\label{B(w, z), B(z, w) in A_0}
        B(w, z), B(z, w) \in A_0
\end{equation}
for every $w \in A^N$ and $z \in A_0^N$.

        As before, $A_0^N$ is a subgroup of $A^N$ as an abelian group
with respect to coordinatewise addition, and we can identify the
quotient group $A^N / A_0^N$ with the group $(A / A_0)^N$ of
$N$-tuples of elements of the quotient $A / A_0$, where addition in
$(A / A_0)^N$ is also defined coordinatewise, as usual.  If (\ref{B(w,
z), B(z, w) in A_0}) holds for every $w \in A^N$ and $z \in A_0^N$,
then $B$ induces a mapping
\begin{equation}
        \widetilde{B} : (A / A_0)^N \times (A / A_0)^N \to A / A_0
\end{equation}
which is additive in each variable.  This leads to a group structure
on $(A / A_0)^N \times (A / A_0)$ associated to $\widetilde{B}$, as in
(\ref{(w, s) diamond (z, t) = (w + z, s + t + B(w, z))}).  The
quotient mapping
\begin{equation}
\label{A^N times A to (A / A_0)^N times (A / A_0)}
        A^N \times A \to (A / A_0)^N \times (A / A_0)
\end{equation}
obtained by applying the natural quotient mapping from $A$ onto $A /
A_0$ in each variable is a homomorphism with respect to (\ref{(w, s)
diamond (z, t) = (w + z, s + t + B(w, z))}) on $A^N \times A$ and the
analogue of (\ref{(w, s) diamond (z, t) = (w + z, s + t + B(w, z))})
associated to $\widetilde{B}$ on $(A / A_0)^N \times (A / A_0)$.  Of
course, the kernel of (\ref{A^N times A to (A / A_0)^N times (A /
A_0)}) is $A_0^N \times A_0$.

        Alternatively, suppose that $A_1$ is another abelian group,
and that $\phi_1$ is a homomorphism from $A$ into $A_1$.  Let $\Phi_1$
be the mapping from $A^N$ into $A_1^N$ obtained by applying $\phi_1$
in each coordinate, so that $\Phi_1$ is a homomorphism as a mapping
from $A^N$ into $A_1^N$ as abelian groups with respect to
coordinatewise addition.  Also let $B_1$ be a mapping from $A_1^N
\times A_1^N$ into $A_1$ which is additive in each variable, as
before, and let $\diamond_1$ be the group structure defined on $A_1^N
\times A_1$ associated to $B_1$ as in (\ref{(w, s) diamond (z, t) = (w
+ z, s + t + B(w, z))}).  Consider the mapping $\Phi$ from $A^N \times
A$ into $A_1^N \times A_1$ defined by
\begin{equation}
        \Phi((z, t)) = (\Phi_1(z), \phi_1(t)).
\end{equation}
If
\begin{equation}
\label{B_1(Phi_1(w), Phi_1(z)) = phi_1(B(w, z))}
        B_1(\Phi_1(w), \Phi_1(z)) = \phi_1(B(w, z))
\end{equation}
for every $w, z \in A^N$, then
\begin{eqnarray}
\label{Phi((w, s) diamond (z, t)) = ... = Phi((w, s)) diamond_1 Phi((z, t))}
 \Phi((w, s) \diamond (z, t)) & = & \Phi((w + z, s + t + B(w, z))) \nonumber \\
 & = & (\Phi_1(w) + \Phi_1(z), \phi_1(s) + \phi_1(t) + \phi_1(B(w, z)))
                                                               \nonumber \\
 & = & (\Phi_1(w) + \Phi_1(z), \phi_1(s) + \phi_1(t)
                                          + B_1(\Phi_1(w), \Phi_1(z)))  \\
 & = & \Phi((w, s)) \diamond_1 \Phi((z, t))        \nonumber
\end{eqnarray}
for every $(w, s), (z, t) \in A^N \times A$, which is to say that
$\Phi$ is a homomorphism from $A^N \times A$ equipped with $\diamond$
into $A_1^N \times A_1$ equipped with $\diamond_1$.  If $A_0$ is the
kernel of $\phi_1$, then $A_0^N$ is the kernel of $\Phi_1$, and $A_0^N
\times A_0$ is the kernel of $\Phi$.  It is easy to see that
(\ref{B(w, z), B(z, w) in A_0}) holds under these conditions, because
\begin{equation}
\label{B_1(Phi_1(w), Phi_1(z)) = 0}
        B_1(\Phi_1(w), \Phi_1(z)) = 0
\end{equation}
as soon as $\Phi_1(w) = 0$ or $\Phi_1(z) = 0$.

\section{Rings and subrings}
\label{rings, subrings}
\setcounter{equation}{0}

        Let $R$ be a commutative ring, not necessarily with a
multiplicative identity element.  If $N$ is a positive integer, then
the set $R^N$ of $N$-tuples of elements of $R$ is also a commutative
ring, with respect to coordinatewise addition and multiplication.  As
before, we shall be more concerned with $R^N$ as a module over $R$,
with respect to coordinatewise addition and scalar multiplication.
Let $B(w, z)$ be an $R$-valued function of $w, z \in R^N$, which is
linear in each variable separately, in the sense that
\begin{equation}
        B(w + w', z) = B(w, z) + B(w', z)
\end{equation}
and
\begin{equation}
        B(w, z + z') = B(w, z) + B(w, z')
\end{equation}
for every $w, w', z, z' \in R^N$, and
\begin{equation}
        B(r \, w, z) = B(w, r \, z) = r \, B(w, z)
\end{equation}
for every $w, z \in R^N$ and $r \in R$.  Under these conditions,
\begin{equation}
\label{(w, s) diamond (z, t) = (w + z, s + t + B(w, z)), rings}
        (w, s) \diamond (z, t) = (w + z, s + t + B(w, z))
\end{equation}
defines a group structure on $R^N \times R$, with $(0, 0)$ as the
identity element, as in the previous sections.

        If $\{b_{j, l}\}_{j, l = 1}^N$ is an $N \times N$ matrix with
entries in $R$, then
\begin{equation}
\label{B(w, z) = sum_{j = 1}^N sum_{l = 1}^N b_{j, l} w_j z_l}
        B(w, z) = \sum_{j = 1}^N \sum_{l = 1}^N b_{j, l} \, w_j \, z_l
\end{equation}
is linear in each variable separately, as in the preceding paragraph.
Note that the bilinear mapping discussed in Section \ref{heisenberg
groups} is of this form, with $N = 2 \, n$.  Conversely, if $R$ has a
multiplicative identity element $e$, then every bilinear mapping $B(w,
z)$ from $R^N \times R^N$ into $R$ is of this form.  To see this, let
$e_j$ be the element of $R^N$ with $j$th component equal to $e$, and
all other components equal to $0$.  In this case, $B(w, z)$ can be
expressed as in (\ref{B(w, z) = sum_{j = 1}^N sum_{l = 1}^N b_{j, l}
w_j z_l}), with
\begin{equation}
\label{b_{j, l} = B(e_j, e_l)}
        b_{j, l} = B(e_j, e_l).
\end{equation}

        If $B(w, z)$ is as in (\ref{B(w, z) = sum_{j = 1}^N sum_{l =
1}^N b_{j, l} w_j z_l}), and $R_0$ is a subring of $R$, then
\begin{equation}
\label{B(w, z) in R_0}
         B(w, z) \in R_0
\end{equation}
for every $w, z \in R_0^N$.  This implies that
\begin{equation}
\label{(w, s) diamond (z, t) in R_0^N times R_0}
        (w, s) \diamond (z, t) \in R_0^N \times R_0
\end{equation}
for every $(w, s), (z, t) \in R_0^N \times R_0$, and that
\begin{equation}
\label{(w, s)^{-1} = (-w, -s + B(w, w)) in R_0^N times R_0}
        (w, s)^{-1} = (-w, -s + B(w, w)) \in R_0^N \times R_0
\end{equation}
for every $(w, s) \in R_0^N \times R_0$, so that $R_0^N \times R_0$ is
a subgroup of $R^N \times R$ with respect to (\ref{(w, s) diamond (z,
t) = (w + z, s + t + B(w, z)), rings}), as in the previous section.
Similarly, if $R_0$ is an ideal in $R$, and $B(w, z)$ is as in
(\ref{B(w, z) = sum_{j = 1}^N sum_{l = 1}^N b_{j, l} w_j z_l}), then
\begin{equation}
\label{B(w, z), B(z, w) in R_0}
        B(w, z), B(z, w) \in R_0
\end{equation}
for every $w \in R_0$ and $z \in R$, and hence $R_0^N \times R_0$ is a
normal subgroup of $R^N \times R$ with respect to (\ref{(w, s) diamond
(z, t) = (w + z, s + t + B(w, z)), rings}), as in the previous section.

\section{Dilations}
\label{dilations}
\setcounter{equation}{0}

        Let $R$ be a commutative ring, let $N$ be a positive integer,
and let $B(w, z)$ be an $R$-valued function of $w, z \in R^N$ which is
linear in each variable separately, as in the previous section.  Thus
(\ref{(w, s) diamond (z, t) = (w + z, s + t + B(w, z)), rings})
defines a group on $R^N \times R$, as before.

        If $r \in R$ and $(z, t) \in R^N \times R$, then put
\begin{equation}
\label{delta_r((z, t)) = (r z, r^2 t)}
        \delta_r((z, t)) = (r \, z, r^2 \, t).
\end{equation}
This defines a mapping from $R^N \times R$ into itself, which satisfies
\begin{eqnarray}
\label{delta_r((w, s) diamond (z, t)) = ...}
\delta_r((w, s) \diamond (z, t)) & = & \delta_r((w + z, s + t + B(w, z)) \\
 & = & (r \, w + r \, z, r^2 \, s + r^2 \, t + r^2 \, B(w, z)) \nonumber \\
 & = & (r \, w + r \, z, r^2 \, s + r^2 \, t + B(r \, w, r \, z)) \nonumber \\
 & = & \delta_r((w, s)) \diamond \delta_r((z, t)) \nonumber
\end{eqnarray}
for every $(w, s), (z, t) \in R^N \times R$.  Therefore $\delta_r$ is
a homomorphism on $R^N \times R$ with respect to (\ref{(w, s) diamond
(z, t) = (w + z, s + t + B(w, z)), rings}).  It is also easy to see
that
\begin{equation}
\label{delta_r circ delta_{r'} = delta_{r r'}}
        \delta_r \circ \delta_{r'} = \delta_{r \, r'}
\end{equation}
for every $r, r' \in R$.  

        If there is a multiplicative identity element $e$ in $R$, then
$\delta_e$ is the identity mapping on $R^N \times R$.  If $r \in R$
has a multiplicative inverse $r^{-1}$ in $R$, then it follows that
$\delta_r$ is an invertible mapping on $R^N \times R$, with
$\delta_r^{-1} = \delta_{r^{-1}}$.

         If $R_0$ is a subring of $R$, then
\begin{equation}
        \delta_r((z, t)) \in R_0^N \times R_0
\end{equation}
for every $r \in R_0$ and $(z, t) \in R_0^N \times R_0$.  If $R_0$ is
an ideal in $R$, then this holds for every $r \in R$ and $(z, t) \in
R_0^N \times R_0$, in addition to $r \in R_0$.

        More generally, let $A$ be an abelian group, with the group
operation written additively, and let $\Sigma$ be a commutative
semigroup, with the group operation written multiplicatively.  
Suppose that $\Sigma$ acts on $A$ in such a way that
\begin{equation}
        a \mapsto r \, a
\end{equation}
is a homomorphism on $A$ for each $r \in \Sigma$, and
\begin{equation}
        r \, (r' \, a) = (r \, r') \, a
\end{equation}
for every $r, r' \in \Sigma$ and $a \in A$.  For example, this happens
when $A$ is a module over a commutative ring $R$, with $\Sigma$ equal
to $R$ as a semigroup with respect to multiplication.  In particular,
any abelian group $A$ may be considered as a module over the ring of
integers, where $n \, a$ is the sum of $n$ $a$'s in $A$ for each $a
\in A$ and $n \in {\bf Z}$.  Conversely, if a semigroup $\Sigma$ acts
on an abelian group $A$ in this way, then the action of $\Sigma$ on
$A$ can be extended to the semigroup ring of formal sums of finitely
many element of $\Sigma$, so that $A$ becomes a module over the
semigroup ring associated to $\Sigma$.

        If $\Sigma$ acts on $A$ as in the previous paragraph, then
$\Sigma$ also acts on $A^N$ as an abelian group with respect to
coordinatewise addition, where elements of $\Sigma$ act on elements of
$A^N$ coordinatewise too.  Let $B(w, z)$ be an $A$-valued function of
$w, z \in A^N$ which is additive in each variable separately, as
before.  Suppose also that $B(w, z)$ is linear with respect to the
action of $\Sigma$ on $A$ in each variable separately, in the sense that
\begin{equation}
\label{B(r w, z) = B(w, r z) = r B(w, z), dilations}
        B(r \, w, z) = B(w, r \, z) = r \, B(w, z)
\end{equation}
for each $r \in \Sigma$ and $w, z \in A^N$.  If $\delta_r$ is defined
on $A^N \times A$ as in (\ref{delta_r((z, t)) = (r z, r^2 t)}) for
each $r \in \Sigma$, then the analogues of (\ref{delta_r((w, s)
diamond (z, t)) = ...}) and (\ref{delta_r circ delta_{r'} = delta_{r
r'}}) still hold in this this situation.  Thus $\delta_r$ defines an
action of $\Sigma$ on $A^N \times A$ as a group with respect to the
operation $\diamond$ associated to $B(w, z)$.

        If $\Sigma$ has an identity element $e$, then one may ask that
$e$ act by the identity mapping on $A$, which is to say that $e \, a =
a$ for each $a \in A$.  In this case, $\delta_e$ also acts by the
identity mapping on $A^N \times A$.

\section{Another subgroup}
\label{another subgroup}
\setcounter{equation}{0}

        Let $A$ be an abelian group, let $N$ be a positive integer,
and let $B(w, z)$ be an $A$-valued function on defined for $w, z \in
A^N$ that is additive in each variable separately, as before.  Thus
$A^N \times A$ is a group with respect to (\ref{(w, s) diamond (z, t)
= (w + z, s + t + B(w, z))}), as usual.  It is easy to see that $\{0\}
\times A$ is a subgroup of $A^N \times A$ with respect to (\ref{(w, s)
diamond (z, t) = (w + z, s + t + B(w, z))}), and that
\begin{equation}
\label{t mapsto (0, t)}
        t \mapsto (0, t)
\end{equation}
defines as isomorphism from $A$ onto $\{0\} \times A$ as a subgroup of
$A^N \times A$ with respect to (\ref{(w, s) diamond (z, t) = (w + z, s
+ t + B(w, z))}).  If $w \in A^N$, then we have seen that
\begin{equation}
\label{B(w, 0) = B(0, w) = 0}
        B(w, 0) = B(0, w) = 0,
\end{equation}
and hence
\begin{equation}
\label{(w, s) diamond (0, t) = (0, t) diamond (w, s) = (w, s + t)}
        (w, s) \diamond (0, t) = (0, t) \diamond (w, s) = (w, s + t)
\end{equation}
for every $s, t \in A$.  This shows that every element of $A^N \times
A$ commutes with every element of $\{0\} \times A$ with respect to
(\ref{(w, s) diamond (z, t) = (w + z, s + t + B(w, z))}).

        In particular, $\{0\} \times A$ is a normal subgroup of $A^N
\times A$ with respect to (\ref{(w, s) diamond (z, t) = (w + z, s + t
+ B(w, z))}).  Consider the natural projection mapping $\Phi$ from
$A^N \times A$ onto $A^N$, defined by
\begin{equation}
        \Phi((z, t)) = z.
\end{equation}
It is easy to see that $\Phi$ is a homomorphism from $A^N \times A$
equipped with (\ref{(w, s) diamond (z, t) = (w + z, s + t + B(w, z))})
onto $A^N$ as an abelian group with coordinatewise addition.  Of
course, the kernel of $\Phi$ is equal to $\{0\} \times A$.

        Suppose now that $R$ is a commutative ring, and that $B(w, z)$
is an $R$-valued function of $w, z \in R^N$ that is linear in each
variable separately, as in the previous sections.  Let $\delta_r$ be
the ``parabolic'' dilation mapping on $R^N \times R$ associated to $r
\in R$ as in the preceding section.  Thus
\begin{equation}
        \delta_r((0, t)) = (0, r^2 \, t)
\end{equation}
for every $r, t \in R$, so that
\begin{equation}
        \delta_r(\{0\} \times R) \subseteq \{0\} \times R
\end{equation}
for every $r \in R$.  We also have that
\begin{equation}
\label{Phi(delta_r((z, t))) = Phi((r z, r^2 t)) = r z}
        \Phi(\delta_r((z, t))) = \Phi((r \, z, r^2 \, t)) = r \, z
\end{equation}
for every $z \in R^N$ and $r, t \in R$, where $r \, z$ corresponds to
coordinatewise scalar multiplication on $R^N$.  There are analogous
statements for the more general situation of a commutative semigroup
$\Sigma$ acting on an abelian group $A$, as in the previous section.

\section{Cocycles}
\label{cocycles}
\setcounter{equation}{0}

        Let $A$ be an abelian group, and let $N$ be a positive
integer, so that $A^N$ is also an abelian group with respect to
coordinatewise addition, as usual.  An $A$-valued function $B(w, z)$
of $w, z \in A^N$ is said to be a $2$-cocycle or factor set if
\begin{equation}
\label{B(z, z') + B(z + z', z'') = B(z, z' + z'') + B(z', z'')}
        B(z, z') + B(z + z', z'') = B(z, z' + z'') + B(z', z'')
\end{equation}
for every $z, z', z'' \in A^N$.  As in Section \ref{heisenberg
groups}, this is exactly the condition that ensures that
\begin{equation}
\label{(w, s) diamond (z, t) = (w + z, s + t + B(w, z)), cocycles}
        (w, s) \diamond (z, t) = (w + z, s + t + B(w, z))
\end{equation}
satisfies the associative law on $A^N \times A$.  In particular, this
condition is satisfied when $B(w, z)$ is an additive function of each
variable separately, as before.

        Note that the sum of two $2$-cocycles is a $2$-cocycle.  If
$B(w, z) = 0$ for every $w, z \in A^N$, then $B(w, z)$ is obviously a
$2$-cocycle.  If $B(w, z)$ is a $2$-cocycle, then $- B(w, z)$ is too,
where the minus sign refers to additive inverses in $A$.  Thus the
collection of $2$-cocycles is a group with respect to addition.

        An $A$-valued function $B(w, z)$ of $w, z \in A^N$ is said to
be a $2$-coboundary if there is an $A$-valued function $b(z)$ of $z
\in A^N$ such that
\begin{equation}
        B(w, z) = b(w) - b(w + z) + b(z)
\end{equation}
for every $w, z \in A$.  In this case, $B(w, z)$ is also a
$2$-cocycle, because
\begin{eqnarray}
\lefteqn{\quad B(z, z') + B(z + z', z'')} \\
 & = & (b(z) - b(z + z') + b(z')) + (b(z + z') - b(z + z' + z'') + b(z''))
                                                       \nonumber \\
 & = & b(z) + b(z') + b(z'') - b(z + z' + z'')         \nonumber \\
 & = & (b(z) - b(z + z' + z'') + b(z' + z'')) + (b(z') - b(z' + z'') + b(z''))
                                                       \nonumber \\
 & = & B(z, z' + z'') + B(z', z'')                     \nonumber
\end{eqnarray}
for every $z, z', z'' \in A^N$.  If $B(w, z) = 0$ for every $w, z \in
A^N$, then $B(w, z)$ is a $2$-coboundary, corresponding to $b(z) = 0$
for each $z \in A^N$.  The sum of two $2$-coboundaries is a
$2$-coboundary, and $- B(w, z)$ is a $2$-coboundary when $B(w, z)$ is
a $2$-coboundary, so that the set of $2$-counbaries is a subgroup of
the group of $2$-cocycles.  The quotient of the group of $2$-cocycles
by the subgroup of $2$-coboundaries is the $2$-dimensional cohomology
group $H^2(A^N, A)$.

        Let $B(w, z)$ be a $2$-cocycle, let $b_0(z)$ be an $A$-valued
function of $z \in A^N$, and let
\begin{equation}
\label{B_0(w, z) = b_0(w) - b_0(w + z) + b_0(z)}
        B_0(w, z) = b_0(w) - b_0(w + z) + b_0(z)
\end{equation}
be the $2$-coboundary associated to $b_0(z)$.  Thus
\begin{equation}
\label{widetilde{B}(w, z) = B(w, z) + B_0(w, z)}
        \widetilde{B}(w, z) = B(w, z) + B_0(w, z)
\end{equation}
is a $2$-cocycle that corresponds to the same element of $H^2(A^N, A)$
as $B(w, z)$.  Let $\widetilde{\diamond}$ be the binary operation on
$A^N \times A$ that corresponds to $\widetilde{B}(w, z)$ as in
(\ref{(w, s) diamond (z, t) = (w + z, s + t + B(w, z)), cocycles}),
and put
\begin{equation}
\label{psi_0((z, t)) = (z, t - b_0(z))}
        \psi_0((z, t)) = (z, t - b_0(z))
\end{equation}
for each $(z, t) \in A^N \times A$.  This defines a one-to-one mapping
from $A^N \times A$ onto itself, which satisfies
\begin{eqnarray}
\label{psi_0((w, s) diamond (z, t)) = ...}
 \psi_0((w, s) \diamond (z, t)) & = & \psi_0((w + z, s + t + B(w, z))) \\
                & = & (w + z, s + t + B(w, z) - b_0(w + z)) \nonumber \\
 & = & (w + z, s + t + \widetilde{B}(w, z) - b_0(w) - b_0(z)) \nonumber \\
 & = & \psi_0((w, s)) \widetilde{\diamond} \psi_0((z, t)) \nonumber
\end{eqnarray}
for every $(w, s), (z, t) \in A^n \times A$.  This corresponds to the
discussion in Section \ref{changing variables} when $b_0(z) = - C(z,
z)$ for an $A$-valued function $C(w, z)$ of $w, z \in A^N$ which is
additive in each variable separately.

        Let $B(w, z)$ be a $2$-cocycle again, and observe that
\begin{equation}
        B(0, 0) + B(0, z'') = B(0, z'') + B(0, z'')
\end{equation}
for every $z'' \in A^N$, by taking $z = z' = 0$ in (\ref{B(z, z') +
B(z + z', z'') = B(z, z' + z'') + B(z', z'')}).  This implies that
\begin{equation}
\label{B(0, z'') = B(0, 0)}
        B(0, z'') = B(0, 0)
\end{equation}
for every $z'' \in A^N$, and similarly we get that
\begin{equation}
\label{B(z, 0) = B(0, 0)}
        B(z, 0) = B(0, 0)
\end{equation}
for every $z \in A^N$, by taking $z' = z'' = 0$ in (\ref{B(z, z') +
B(z + z', z'') = B(z, z' + z'') + B(z', z'')}).  It follows that
\begin{equation}
 (w, s) \diamond (0, - B(0, 0)) = (0, - B(0, 0)) \diamond (w, s) = (w, s)
\end{equation}
for every $(w, s) \in A^N \times A$, so that $(0, - B(0, 0))$ is the
identity element for $\diamond$.

        If $z = -z' = z''$, then (\ref{B(z, z') + B(z + z', z'') =
B(z, z' + z'') + B(z', z'')}) reduces to
\begin{equation}
        B(z, -z) + B(0, z) = B(z, 0) + B(-z, z).
\end{equation}
This implies that
\begin{equation}
\label{B(z, -z) = B(-z, z)}
        B(z, -z) = B(-z, z)
\end{equation}
for every $z \in A^N$, because of (\ref{B(0, z'') = B(0, 0)}) and
(\ref{B(z, 0) = B(0, 0)}).  Hence
\begin{eqnarray}
\label{(w, s) diamond (-w, -s - B(w, -w) - B(0, 0)) = ... = (0, - B(0, 0))}
\lefteqn{(w, s) \diamond (-w, -s - B(w, -w) - B(0, 0))} \\
           & = & (-w, -s - B(w, -w) - B(0, 0)) \diamond (w, s) \nonumber \\
           & = & (0, - B(0, 0)) \nonumber
\end{eqnarray}
for every $(w, s) \in A^N \times A$, which shows that $(-w, -s - B(w,
-w) - B(0, 0))$ is the inverse of $(w, s)$ with respect to $\diamond$.

        Of course, these formulae would be a bit nicer if $B(0, 0) =
0$.  If $b_0(z)$ is a constant function on $A^N$, then the
corresponding $2$-coboundary is a constant function on $A^N \times
A^N$, with the same constant value.  This implies that
\begin{equation}
\label{widehat{B}(w, z) = B(w, z) - B(0, 0)}
        \widehat{B}(w, z) = B(w, z) - B(0, 0)
\end{equation}
is a $2$-cocycle that is cohomologous to $B(w, z)$ and satisfies
$\widehat{B}(0, 0) = 0$.

        At any rate, if $B(w, z)$ is a $2$-cocycle, then then it
follows that $A^N \times A$ is a group with respect to $\diamond$.  As
in the previous section, $\{0\} \times A$ is a subgroup of $A^N \times
A$ with respect to $\diamond$, and
\begin{equation}
\label{t mapsto (0, t - B(0, 0))}
        t \mapsto (0, t - B(0, 0))
\end{equation}
defines an isomorphism from $A$ with its original group structure onto
$\{0\} \times A$ as a subgroup of $A^N \times A$ equipped with
$\diamond$.  We also have that
\begin{equation}
\label{(w, s) diamond (0, t) = (0, t) diamond (w, s) = (w, s + t + B(0, 0))}
 (w, s) \diamond (0, t) = (0, t) \diamond (w, s) = (w, s + t + B(0, 0))
\end{equation}
for every $w \in A^N$ and $s, t \in A$, by (\ref{B(0, z'') = B(0, 0)})
and (\ref{B(z, 0) = B(0, 0)}).  In particular, $\{0\} \times A^N$ is a
normal subgroup of $A^N \times A$ with respect to $\diamond$.  If
$\Phi$ is the natural projection from $A^N \times A$ onto $A^N$, then
$\Phi$ is a homomorphism from $A^N \times A$ equipped with $\diamond$
onto $A^N$ as an abelian group with respect to coordinatewise
addition, and $\{0\} \times A$ is the kernel of $\Phi$, as before.

\section{Haar measure on ${\bf Q}_p$}
\label{haar measure on Q_p}
\setcounter{equation}{0}

        If $p$ is a prime number, then the field of $p$-adic numbers
${\bf Q}_p$ is actually a topological field with respect to the
topology determined by the $p$-adic metric, in the sense that the
field operations are continuous.  This is easy to check, in the same
way as for the real numbers.  In particular, the $p$-adic numbers form
a locally compact abelian group with respect to addition.  This
implies that there is a positive Borel measure $\mu$ on ${\bf Q}_p$
such that $\mu(K) < \infty$ when $K \subseteq {\bf Q}_p$ is compact,
$\mu(U) > 0$ when $U \subseteq {\bf Q}_p$ is nonempty and open, and
$\mu$ is invariant under translations.  More precisely, if $E
\subseteq {\bf Q}_p$ and $y \in {\bf Q}_p$, then we put
\begin{equation}
        E + y = \{x + y : x \in E\}.
\end{equation}
If $E$ is a Borel set in ${\bf Q}_p$, then $E + y$ is also a Borel set
for each $y \in {\bf Q}_p$, because
\begin{equation}
        x \mapsto x + y
\end{equation}
is a homeomorphism on ${\bf Q}_p$, by continuity of addition on ${\bf
Q}_p$.  The condition that $\mu$ be invariant under translations means that
\begin{equation}
        \mu(E + y) = \mu(E)
\end{equation}
for every Borel set $E \subseteq {\bf Q}_p$ and $y \in {\bf Q}_p$.

        Such a measure $\mu$ on ${\bf Q}_p$ is known as Haar measure,
and it is unique up to multiplication by a positive real number.
Normally for a locally compact abelian topological group, one should
also ask $\mu$ to satisfy some additional regularity conditions, to
get this uniqueness statement in particular.  In this case, as on the
real line, it is well known that positive Borel measures that are
finite on compact sets automatically satisfy some regularity
conditions, because every open set is $\sigma$-compact.  The latter
can be derived from the fact that the topology is determined by a
metric, and that closed and bounded sets with respect to the metric
are compact.

        To get uniqueness of Haar measure $\mu$ on ${\bf Q}_p$, let us
also ask that
\begin{equation}
        \mu({\bf Z}_p) = 1.
\end{equation}
This implies that
\begin{equation}
\label{mu(p^j {bf Z}_p) = p^{-j}}
        \mu(p^j \, {\bf Z}_p) = p^{-j}
\end{equation}
for each positive integer $j$, because ${\bf Z}_p$ can be expressed as
the union of $p^j$ pairwise disjoint translates of $p^j \, {\bf Z}_p$.
Similarly, if $j$ is a negative integer, then $p^j \, {\bf Z}_p$ can
be expressed as the union of $p^{-j}$ pairwise disjoint translates of
${\bf Z}_p$, so that (\ref{mu(p^j {bf Z}_p) = p^{-j}}) holds for all
integers $j$.  It follows that
\begin{equation}
\label{mu(p^j {bf Z}_p + y) = p^{-j}}
        \mu(p^j \, {\bf Z}_p + y) = p^{-j}
\end{equation}
for each $j \in {\bf Z}$ and $y \in {\bf Q}_p$, because $\mu$ is
translation-invariant.  This basically characterizes Haar measure on
${\bf Q}_p$, in the same way that Lebesgue measure on the real line is
characterized by the length of an interval.

        Alternatively, one can look at the corresponding Haar integral
\begin{equation}
\label{I(f) = int_{{bf Q}_p} f d mu}
        I(f) = \int_{{\bf Q}_p} f \, d\mu
\end{equation}
as a linear functional on the vector space $C_{com}({\bf Q}_p)$ of
continuous real-valued functions with compact support on ${\bf Q}_p$,
which is nonnegative in the sense that $I(f) \ge 0$ when $f(x) \ge 0$
for every $x \in {\bf Q}_p$.  More precisely,
\begin{equation}
        I(f) > 0
\end{equation}
when $f(x) \ge 0$ for every $x \in {\bf Q}_p$ and $f(x) > 0$ for some
$x \in {\bf Q}_p$, by the positivity properties of Haar measure, and
\begin{equation}
        I(f) = I(g)
\end{equation}
when $g(x) = f(x + y)$ for some $y \in {\bf Q}_p$ and every $x \in
{\bf Q}_p$, by the translation-invariance of Haar measure.
Conversely, if $I(f)$ is a nonnegative translation-invariant linear
functional on $C_{com}({\bf Q}_p)$, then the Riesz representation
theorem leads to a translation-invariant nonnegative Borel measure
$\mu$ on ${\bf Q}_p$ such that (\ref{I(f) = int_{{bf Q}_p} f d mu})
holds for every $f \in C_{com}({\bf Q}_p)$.  Thus one way to get Haar
measure on ${\bf Q}_p$ is to first define $I(f)$ on $C_{com}({\bf
Q}_p)$ with these properties, using Riemann sums and (\ref{mu(p^j {bf
Z}_p + y) = p^{-j}}), and then apply the Riesz representation theorem.
This is analogous to obtaining Lebesgue measure on the real line from
the ordinary Riemann integral as a nonnegative linear functional on
the vector space $C_{com}({\bf R})$ of continuous real-valued
functions with compact support on ${\bf R}$ using the Riesz
representation theorem.

\section{Haar measure on ${\bf Q}_p^n$}
\label{haar measure on Q_p^n}
\setcounter{equation}{0}

        Let $p$ be a prime number, and let $n$ be a positive integer.
As usual, we can consider ${\bf Q}_p^n$ as a vector space over ${\bf
Q}_p$ with respect to coordinatewise addition and scalar
multiplication.  This is also a topological space with respect to the
product topology, using the standard topology on ${\bf Q}_p$
determined by the $p$-adic metric on each factor.  It is easy to see
that the vector space operations are continuous with respect to the
product topology, so that ${\bf Q}_p^n$ is a locally compact abelian
group with respect to coordinatewise addition in particular.  Let
$\mu$ be the normalized Haar measure on ${\bf Q}_p$, and let $\mu_n$
be the corresponding product measure on ${\bf Q}_p^n$, using $\mu$ on
each factor.  This satisfies the analogous positivity and
translation-invariance properties on ${\bf Q}_p^n$, and
\begin{equation}
        \mu_n({\bf Z}_p^n) = \mu({\bf Z}_p)^n = 1.
\end{equation}
Thus $\mu_n$ is normalized Haar measure on ${\bf Q}_p^n$, as a locally
compact abelian group with respect to coordinatewise addition.

        Suppose that $n \ge 2$, and let $\phi(x_2, \ldots, x_n)$ be a
continuous ${\bf Q}_p$-valued function of $x_2, \ldots, x_n \in {\bf
Q}_p$.  Consider
\begin{equation}
 \Phi(x_1, x_2, \ldots, x_n) = (x_1 + \phi(x_2, \ldots, x_n), x_2, \ldots, x_n)
\end{equation}
as a mapping from ${\bf Q}_p^n$ into itself.  More precisely, this is
a one-to-one mapping from ${\bf Q}_p^n$ onto itself, with inverse
given by
\begin{equation}
        \Phi^{-1}(x_1, x_2, \ldots, x_n)
          = (x_1 - \phi(x_2, \ldots, x_n), x_2, \ldots, x_n).
\end{equation}
Note that $\Phi$ and $\Phi^{-1}$ are both continuous as mappings from
${\bf Q}_p^n$ into itself.

        If $f$ is a nonnegative real-valued Borel measurable function
on ${\bf Q}_p^n$, then
\begin{eqnarray}
\lefteqn{\int_{{\bf Q}_p} f(x_1 + \phi(x_2, \ldots, x_n), x_2, \ldots, x_n)
                                          \, d\mu(x_1)} \\
         & = & \int_{{\bf Q}_p} f(x_1, x_2, \ldots, x_n) \, d\mu(x_1) \nonumber
\end{eqnarray}
for all $x_2, \ldots, x_n \in {\bf Q}_p$, by the
translation-invariance of Haar measure on ${\bf Q}_p$.  Integrating
both sides of this equation in $x_2, \ldots, x_n$, we get that
\begin{equation}
\label{int_{{bf Q}_p^n} f circ Phi d mu_n = int_{{bf Q}_p^n} f d mu}
 \int_{{\bf Q}_p^n} f \circ \Phi \, d\mu_n = \int_{{\bf Q}_p^n} f \, d\mu,
\end{equation}
which also holds for integrable real-valued functions on ${\bf Q}_p^n$.
If we take $f$ to be the indicator function associated to a Borel set
$E \subseteq {\bf Q}_p^n$, then we get that
\begin{equation}
        \mu_n(\Phi^{-1}(E)) = \mu_n(E).
\end{equation}
Equivalently,
\begin{equation}
        \mu_n(\Phi(A)) = \mu_n(A)
\end{equation}
for every Borel set $A \subseteq {\bf Q}_p^n$.

\section{Haar measure on $H_n({\bf Q}_p)$}
\label{haar measure on H_n(Q_p)}
\setcounter{equation}{0}

        Let $p$ be a prime number, let $N$ be a positive integer, and
let $B(w, z)$ be a bilinear form on ${\bf Q}_p^N$, which is to say a
${\bf Q}_p$-valued function of $w, z \in {\bf Q}_p^N$ which is linear
in each variable separately.  As usual,
\begin{equation}
\label{(w, s) diamond (z, t) = (w + z, s + t + B(w, z)), haar measure}
        (w, s) \diamond (z, t) = (w + z, s + t + B(w, z))
\end{equation}
defines a group operation on ${\bf Q}_p^N \times {\bf Q}_p$, which includes
the situation in Section \ref{heisenberg groups} with $N = 2 \, n$.
Remember that $(0, 0)$ is the identity element of ${\bf Q}_p^N \times 
{\bf Q}_p$, and that the inverse of $(w, s) \in {\bf Q}_p^N \times {\bf Q}_p$
is given by
\begin{equation}
\label{(w, s)^{-1} = (-w, -s + B(w, w)), haar measure}
        (w, s)^{-1} = (-w, -s + B(w, w)).
\end{equation}
If $B(w, z)$ is antisymmetric in $w$ and $z$, in the sense that
\begin{equation}
\label{B(z, w) = - B(w, z), haar measure}
        B(z, w) = - B(w, z)
\end{equation}
for every $w, z \in {\bf Q}_p^N$, then $B(w, w) = 0$ for every $w \in
{\bf Q}_p^N$, and
\begin{equation}
\label{(w, s)^{-1} = (-w, -s), haar measure}
        (w, s)^{-1} = (-w, -s)
\end{equation}
for every $(w, s) \in {\bf Q}_p^N \times {\bf Q}_p$.

        It is easy to see that the group operations on ${\bf Q}_p^N
\times {\bf Q}_p$ are continuous with respect to the product topology
associated to the standard topology on ${\bf Q}_p$, determined by the
$p$-adic metric.  Thus ${\bf Q}_p^N \times {\bf Q}_p$ is a locally
compact topological group, even if it is not necessary abelian.  Let
$\mu$ be normalized Haar measure on ${\bf Q}_p$, and let $\mu_{N + 1}$
be the corresponding product measure on ${\bf Q}_p^N \times {\bf
Q}_p$, as in the previous sections.  If $E \subseteq {\bf Q}_p^N
\times {\bf Q}_p$ is a Borel set, then
\begin{equation}
        E \diamond (z, t) = \{(w, s) \diamond (z, t) : (w, s) \in E\}
\end{equation}
and
\begin{equation}
        (z, t) \diamond E = \{(z, t) \diamond (w, s) : (w, s) \in E\}
\end{equation}
are also Borel sets in ${\bf Q}_p^N \times {\bf Q}_p$ for every $(z,
t) \in {\bf Q}_p^N \times {\bf Q}_p$.  This follows from the
continuity of the group operation (\ref{(w, s) diamond (z, t) = (w +
z, s + t + B(w, z)), haar measure}), which implies that
\begin{equation}
        (w, s) \mapsto (w, s) \diamond (z, t) \quad\hbox{and}\quad
         (w, s) \mapsto (z, t) \diamond (w, s)
\end{equation}
are homeomorphisms on ${\bf Q}_p^N \times {\bf Q}_p$ for every $(z, t)
\in {\bf Q}_p^N \times {\bf Q}_p$.  Moreover,
\begin{equation}
\label{mu_{N + 1}(E diamond (z, t)) = ... = mu_{N + 1}(E)}
        \mu_{N + 1}(E \diamond (z, t)) = \mu_{N + 1}((z, t) \diamond E)
                                       = \mu_{N + 1}(E)
\end{equation}
for every Borel set $E \subseteq {\bf Q}_p^n \times {\bf Q}_p$ and
$(z, t) \in {\bf Q}_p^N \times {\bf Q}_p$, by the computations
discussed in the previous section.  Thus $\mu_{N + 1}$ satisfies the
requirements of Haar measure on ${\bf Q}_p^N \times {\bf Q}_p$ with
respect to (\ref{(w, s) diamond (z, t) = (w + z, s + t + B(w, z)),
haar measure}), for both translations on the left and on the right.

\section{Haar measure and dilations}
\label{haar measure, dilations}
\setcounter{equation}{0}

        Let $p$ be a prime number, and put
\begin{equation}
        r \, E = \{r \, x : x \in E\}
\end{equation}
for each $r \in {\bf Q}_p$ and $E \subseteq {\bf Q}_p$.  If $E$ is a
Borel set in ${\bf Q}_p$, then $r \, E$ is also a Borel set, because
$r \, E = \{0\}$ when $r = 0$ and $E \ne \emptyset$, and because
\begin{equation}
        x \mapsto r \, x
\end{equation}
is a homeomorphism on ${\bf Q}_p$ when $r \ne 0$.  Let $\mu$ be
normalized Haar measure on ${\bf Q}_p$, and let us check that
\begin{equation}
\label{mu(r E) = |r|_p mu(E)}
        \mu(r \, E) = |r|_p \, \mu(E)
\end{equation}
for every $r \in {\bf Q}_p$ and Borel set $E \subseteq {\bf Q}_p$.
This is trivial when $r = 0$, and otherwise $\mu(r \, E)$ is a Borel
measure on ${\bf Q}_p$ with the same positivity and
translation-invariance properties as Haar measure, which implies that
$\mu(r \, E)$ is equal to a positive constant multiple of $\mu(E)$.
To determine the constant, it suffices to compute $\mu(r \, {\bf
Z}_p)$.  If $r \ne 0$, then $|r|_p = p^{-j}$ for some integer $j$, and
it is easy to see that
\begin{equation}
        r \, {\bf Z}_p = p^j \, {\bf Z}_p.
\end{equation}
Hence $\mu(r \, {\bf Z}_p) = \mu(p^j \, {\bf Z}_p) = p^{-j}$, by
(\ref{mu(p^j {bf Z}_p) = p^{-j}}), which implies (\ref{mu(r E) = |r|_p
mu(E)}), because $\mu({\bf Z}_p) = 1$.

        Similarly, if $r \in {\bf Q}_p$ and $E \subseteq {\bf Q}_p^n$
for some positive integer $n$, then $r \, E$ can be defined as before,
using coordinatewise scalar multiplication on ${\bf Q}_p^n$.  If $E$
is a Borel set in ${\bf Q}_p^n$, then $r \, E$ is also a Borel set in
${\bf Q}_p^n$, for the same reasons as before.  Let $\mu_n$ be the
product measure on ${\bf Q}_p^n$ corresponding to $\mu$ on each
factor.  In this case,
\begin{equation}
\label{mu_n(r E) = |r|_p^n mu_n(E)}
        \mu_n(r \, E) = |r|_p^n \, \mu_n(E).
\end{equation}
This can be shown by an argument like the one in the previous
paragraph, or by observing that $\mu_n(r \, E)$ is the same as the
product of $n$ copies of the corresponding dilation of $\mu$ on ${\bf
Q}_p$.

        Now let $N$ be a positive integer, and put
\begin{equation}
\label{delta_r((z, t)) = (r z, r^2 t), haar measure}
        \delta_r((z, t)) = (r \, z, r^2 \, t)
\end{equation}
for each $r \in {\bf Q}_p$ and $(z, t) \in {\bf Q}_p^N \times {\bf
Q}_p$, as in Section \ref{dilations}.  If $E \subseteq {\bf Q}_p^N
\times {\bf Q}_p$ is a Borel set, then $\delta_r(E)$ is a Borel set in
${\bf Q}_p^N \times {\bf Q}_p$ for each $r \in {\bf Q}_p$, for the
same reasons as before.  One can also check that
\begin{equation}
\label{mu_{N + 1}(delta_r(E)) = |r|_p^{N + 2} mu_{N + 1}(E)}
        \mu_{N + 1}(\delta_r(E)) = |r|_p^{N + 2} \, \mu_{N + 1}(E)
\end{equation}
for every Borel set $E \subseteq {\bf Q}_p^N \times {\bf Q}_p$ and $r
\in {\bf Q}_p$.  Here $\mu_{N + 1}$ is the product measure on ${\bf
Q}_p^N \times {\bf Q}_p$ associated to $\mu$ on each factor, and
(\ref{mu_{N + 1}(delta_r(E)) = |r|_p^{N + 2} mu_{N + 1}(E)}) follows
by identifying $\mu_{N + 1}(\delta_r(E))$ with a product of
appropriate dilations of $\mu$ on ${\bf Q}_p$.

\section{An invariant ultrametric}
\label{invariant ultrametric}
\setcounter{equation}{0}

        Let $p$ be a prime number, let $N$ be a positive integer,
and let
\begin{equation}
\label{B(w, z) = sum_{j = 1}^N sum_{l = 1}^N b_{j, l} w_j z_l, Q_p^N}
        B(w, z) = \sum_{j = 1}^N \sum_{l = 1}^N b_{j, l} w_j \, z_l
\end{equation}
be a bilinear form on ${\bf Q}_p^N$.  In this section, it will be
convenient to ask that
\begin{equation}
\label{b_{j, l} in {bf Z}_p}
        b_{j, l} \in {\bf Z}_p
\end{equation}
for each $j, l = 1, \ldots, N$, which implies that
\begin{equation}
\label{|B(w, z)|_p le (max_{1 le j le n} |w_j|_p) (max_{1 le l le n} |z_l|_p)}
        |B(w, z)|_p \le \Big(\max_{1 \le j \le n} |w_j|_p\Big)
                         \Big(\max_{1 \le l \le n} |z_l|_p\Big)
\end{equation}
for every $w, z \in {\bf Q}_p^N$.  Also let $(w, s) \diamond (z, t)$
be the corresponding group structure on ${\bf Q}_p^N \times {\bf
Q}_p$, as in (\ref{(w, s) diamond (z, t) = (w + z, s + t + B(w, z)),
haar measure}).  If $(z, t) \in {\bf Q}_p^N \times {\bf Q}_p$, then put
\begin{equation}
\label{||(z, t)|| = max(|z_1|_p, ldots, |z_N|_p, |t|_p^{1/2})}
        \|(z, t)\| = \max(|z_1|_p, \ldots, |z_N|_p, |t|_p^{1/2}).
\end{equation}
Thus $\|(z, t)\| \ge 0$ for every $(z, t) \in {\bf Q}_p^N \times {\bf
Q}_p$, and $\|(z, t)\| = 0$ if and only if $(z, t) = (0, 0)$.

        Observe that
\begin{equation}
        |B(w, z)|_p \le \max (|w_1|_p^2, \ldots, |w_n|_p^2,
                                 |z_1|_p^2, \ldots, |z_n|_p^2)
\end{equation}
for every $w, z \in {\bf Q}_p^n$, by (\ref{|B(w, z)|_p le (max_{1 le j
le n} |w_j|_p) (max_{1 le l le n} |z_l|_p)}).  This implies that
\begin{equation}
 \quad |s + t + B(w, z)|_p \le \max (|s|_p, |t|_p,
              |w_1|_p^2, \ldots, |w_n|_p^2, |z_1|_p^2, \ldots, |z_n|_p^2)
\end{equation}
and hence
\begin{equation}
\label{||(w, s) diamond (z, t)|| le max (||(w, s)||, ||(z, t)||)}
        \|(w, s) \diamond (z, t)\| \le \max (\|(w, s)\|, \|(z, t)\|)
\end{equation}
for every $(w, s), (z, t) \in {\bf Q}_p^N \times {\bf Q}_p$.  One can
also check that
\begin{equation}
\label{||(w, s)^{-1}|| = ||(w, s)||}
        \|(w, s)^{-1}\| = \|(w, s)\|
\end{equation}
for every $(w, s) \in {\bf Q}_p^N \times {\bf Q}_p$.  This is very
easy to see when $B(w, z)$ is antisymmetric, so that $B(w, w) = 0$ for
every $w \in {\bf Q}_p^N$, and otherwise one can show this by
estimating $|B(w, w)|_p$ by $\max (|w_1|_p^2, \ldots, |w_n|_p^2)$, as
before.

        If $(w, s), (z, t) \in {\bf Q}_p^N \times {\bf Q}_p$, then put
\begin{equation}
\label{d((w, s), (z, t)) = ||(w, s)^{-1} diamond (z, t)||}
        d((w, s), (z, t)) = \|(w, s)^{-1} \diamond (z, t)\|.
\end{equation}
Of course,
\begin{equation}
        ((w, s)^{-1} \diamond (z, t))^{-1} = (z, t)^{-1} \diamond (w, s),
\end{equation}
so that
\begin{equation}
\label{d((w, s), (z, t)) = d((z, t), (w, s))}
        d((w, s), (z, t)) = d((z, t), (w, s)),
\end{equation}
for every $(w, s), (z, t) \in {\bf Q}_p^N \times {\bf Q}_p$, by
(\ref{||(w, s)^{-1}|| = ||(w, s)||}).  One can also check that
(\ref{d((w, s), (z, t)) = ||(w, s)^{-1} diamond (z, t)||}) satisfies
the ultrametric version of the triangle inequality, using (\ref{||(w,
s) diamond (z, t)|| le max (||(w, s)||, ||(z, t)||)}).  Thus
(\ref{d((w, s), (z, t)) = ||(w, s)^{-1} diamond (z, t)||}) defines an
ultrametric on ${\bf Q}_p^N \times {\bf Q}_p$, which is easily seen to
determine the same topology on ${\bf Q}_p^N \times {\bf Q}_p$ as the
product topology associated to the standard topology on ${\bf Q}_p$.

        By construction, this ultrametric is invariant under left
translations on ${\bf Q}_p^N \times {\bf Q}_p$, in the sense that
\begin{equation}
\label{d((u, v) diamond (w, s), (u, v) diamond (z, t)) = d((w, s), (z, t))}
        d((u, v) \diamond (w, s), (u, v) \diamond (z, t)) = d((w, s), (z, t))
\end{equation}
for every $(u, v), (w, s), (z, t) \in {\bf Q}_p^N \times {\bf Q}_p$.
More precisely, this follows from the algebraic fact that
\begin{eqnarray}
\lefteqn{((u, v) \diamond (w, s))^{-1} \diamond ((u, v) \diamond (z, t))} \\
 & = & (w, s)^{-1} \diamond (u, v)^{-1} \diamond (u, v) \diamond (z, t)
                                                             \nonumber \\
 & = & (w, s)^{-1} \diamond (z, t) \nonumber
\end{eqnarray}
for every $(u, v), (w, s), (z, t) \in {\bf Q}^N \times {\bf Q}_p$.
In addition, note that
\begin{equation}
\label{||delta_r((z, t))|| = |r|_p ||(z, t)||}
        \|\delta_r((z, t))\| = |r|_p \, \|(z, t)\|
\end{equation}
for every $r \in {\bf Q}_p$, where $\delta_r((z, t)) = (r \, z, r^2,
t)$, as in (\ref{delta_r((z, t)) = (r z, r^2 t), haar measure}).
This implies that
\begin{equation}
\label{d(delta_r((w, s)), delta_r((z, t))) = |r|_p d((w, s), (z, t))}
        d(\delta_r((w, s)), \delta_r((z, t))) = |r|_p \, d((w, s), (z, t))
\end{equation}
for every $(w, s), (z, t) \in {\bf Q}_p^N \times {\bf Q}_p$ and $r \in
{\bf Q}_p$.  This also uses the fact that $\delta_r$ is a homomorphism
on ${\bf Q}_p^N \times {\bf Q}_p$ for each $r \in {\bf Q}_p$, as in
Section \ref{dilations}.

\section{Polynomials}
\label{polynomials}
\setcounter{equation}{0}

        Let $R$ be a commutative ring, not necessarily with a
multiplicative identity element.  A polynomial in one variable over
$R$ is a formal expression
\begin{equation}
\label{f(x) = a_l x^l + a_{l - 1} x^{l - 1} + cdots + a_1 x + a_0}
        f(x) = a_l \, x^l + a_{l - 1} \, x^{l - 1} + \cdots + a_1 \, x + a_0,
\end{equation}
where the coefficients $a_j$ are elements of $R$, $j = 0, \ldots, l$.
Here $x$ is an indeterminant, and $x^j$ denotes the formal product of
$j$ $x$'s.  This defines a function on $R$ by substituting particular
elements of $R$ for $x$, and similarly one can get a function on any
ring that contains $R$ as a subring, or on which $R$ acts by
multiplication in a nice way.  One can also take the sum, product, and
composition of two polynomials in the usual way, and this is
compatible with the sum, product, and composition of the corresponding
functions on $R$ or rings on which $R$ acts.

        Now let $n$ be a positive integer, and let $x_1, \ldots, x_n$
be $n$ commuting indeterminants.  If $\alpha$ is a multi-index, which
is to say an $n$-tuple $(\alpha_1, \ldots, \alpha_n)$ of nonnegative
integers, then the corresponding monomial $x^\alpha$ is the formal
product of $x_j^{\alpha_j}$ for $j = 1, \ldots, n$.  Put
\begin{equation}
        |\alpha| = \sum_{j = 1}^n \alpha_j,
\end{equation}
which is the degree of the monomial $x^\alpha$.  A polynomial in $n$
variables over $R$ is a formal expression
\begin{equation}
\label{f(x) = sum_{|alpha| le d} a_alpha x^alpha}
        f(x) = \sum_{|\alpha| \le d} a_\alpha \, x^\alpha,
\end{equation}
where $d$ is a nonnegative integer and $a_\alpha \in R$ for each $\alpha$.
More precisely, the sum extends over the finite set of all multi-indices
$\alpha$ with $|\alpha| \le d$, and the degree of $f$ is less than or
equal to $d$ in this case.  As before, this defines an $R$-valued
function on $R^n$, by substituting specific elements of $R$ for $x_1,
\ldots, x_n$.  If $R_1$ is another commutative ring that contains $R$
as a subring, or on which $R$ acts by multiplication in a nice way,
then one can use $f(x)$ to get an $R_1$-valued function on $R_1^n$ in
the same way.  One can again take sums, products, and compositions of
polynomials, which are compatible with sums, products, and
compositions of the corresponding functions.

        If a polynomial $f(x)$ can be expressed as
\begin{equation}
\label{f(x) = sum_{|alpha| = d} a_alpha x^alpha}
        f(x) = \sum_{|\alpha| = d} a_\alpha \, x^\alpha
\end{equation}
for some nonnegative integer $d$, then we say that $f(x)$ is homogeneous
of degree $d$.  This implies that
\begin{equation}
\label{f(r x) = r^d f(x)}
        f(r \, x) = r^d \, f(x),
\end{equation}
where $r$ is another indeterminant that commutes with $x_1, \ldots,
x_n$.  Of course, (\ref{f(r x) = r^d f(x)}) still holds if we take $r$
to be a particular element of $R$.

        Suppose that $n \ge 2$, and put $N = n - 1$ and
\begin{equation}
        |\alpha|_\delta = \sum_{j = 1}^N \alpha_j + 2 \, \alpha_{N + 1}
\end{equation}
for each multi-index $\alpha$.  If a polynomial $f(x)$ can be expressed as
\begin{equation}
        f(x) = \sum_{|\alpha|_\delta = d} a_\alpha \, x^\alpha
\end{equation}
for some nonnegative integer $d$, then we say that $f(x)$ is
$\delta$-homogeneous of degree $d$.  This implies that
\begin{equation}
\label{f(r x_1, ..., r x_N, r^2 x_{N + 1}) = r^d f(x_1, ..., x_N, x_{N + 1})}
        f(r \, x_1, \ldots, r \, x_N, r^2 \, x_{N + 1})
           = r^d \, f(x_1, \ldots, x_N, x_{N + 1}),
\end{equation}
where $r$ is another indeterminant that commutes with the $x_j$'s, as before.
This is the homogeneity condition that corresponds to the dilations $\delta_r$
on $R^N \times R$ as in Section \ref{dilations}.

\section{Formal power series}
\label{formal power series}
\setcounter{equation}{0}

        Let $R$ be a commutative ring again, not necessarily with a
multiplicative identity element.  A formal power series in one
variable over $R$ is a formal expression
\begin{equation}
\label{f(x) = sum_{j = 0}^infty a_j x^j}
        f(x) = \sum_{j = 0}^\infty a_j \, x^j,
\end{equation}
where $a_j \in R$ for each $j \ge 0$ and $x$ is an indeterminant.  Similarly,
a formal power series in $n$ commuting indeterminants $x_1, \ldots, x_n$
is given by
\begin{equation}
\label{f(x) = sum_alpha a_alpha x^alpha}
        f(x) = \sum_\alpha a_\alpha \, x^\alpha,
\end{equation}
where the sum is taken over all multi-indices $\alpha$, and $a_\alpha
\in R$ for each multi-index $\alpha$.  One can add and multiply formal
power series in the usual way, but one should be more careful about
compositions, to avoid infinite sums of multiples of the same monomial
$x^\alpha$ with nonzero coefficients.

        If $f(x)$ is a formal power series in one variable as in
(\ref{f(x) = sum_{j = 0}^infty a_j x^j}), then the formal derivative
of $f(x)$ is given by
\begin{equation}
\label{f'(x) = sum_{j = 1}^infty j a_j x^{j - 1}}
        f'(x) = \sum_{j = 1}^\infty j \, a_j \, x^{j - 1}.
\end{equation}
More precisely, $j \, a_j$ is the sum of $j$ $a_j$'s in $R$.
Using the binomial theorem, one can check that
\begin{equation}
\label{f(x + h) = f(x) + f'(x) h + g(x, h) h^2}
        f(x + h) = f(x) + f'(x) \, h + g(x, h) \, h^2,
\end{equation}
where $h$ is another indeterminant that commutes with $x$, and $g(x,
h)$ is a formal power series in $x$ and $h$.  If $f(x)$ is a
polynomial, then $f'(x)$ is also a polynomial, and one can take $g(x,
h)$ to be a polynomial in $x$ and $h$.

        Similarly, if $f(x)$ is a formal power series in $n$ commuting
indeterminants $x_1, \ldots, x_n$ as in (\ref{f(x) = sum_alpha a_alpha
x^alpha}), then the formal partial derivative of $f(x)$ in $x_l$ is
defined for $l = 1, \ldots, n$ by
\begin{equation}
\label{partial_l f(x) = frac{partial f}{partial x_l} (x) = ...}
        \partial_l f(x) = \frac{\partial f}{\partial x_l} (x)
          = \sum_{\alpha_l \ge 1} \alpha_l \, a_\alpha \, x^{\alpha'(l)},
\end{equation}
where the sum is taken over all multi-indices $\alpha$ such that
$\alpha_l \ge 1$, and where $\alpha'(l)$ is the multi-index given by
\begin{equation}
        \alpha'_l(l) = \alpha_l - 1, \quad
         \alpha'_k(l) = \alpha_k \hbox{ when } k \ne l.
\end{equation}
If $h_1, \ldots, h_n$ are commuting indeterminants that commute with
$x_1, \ldots, x_n$, then one can use the binomial theorem again to get that
\begin{equation}
\label{f(x + h) = f(x) + sum_{l = 1}^n partial_l f(x) h_l + ...}
        f(x + h) = f(x) + \sum_{l = 1}^n \partial_l f(x) \, h_l
                         + \sum_{k, l = 1}^n g_{k, l}(x, h) \, h_k \, h_l,
\end{equation}
where $g_{k, l}(x, h)$ is a formal power series in $x_1, \ldots, x_n$
and $h_1, \ldots, h_n$ for each $k$, $l$.  If $f(x)$ is a polynomial
in $x_1, \ldots, x_n$, then $\partial f/\partial x_l$ is also a
polynomial for $l = 1, \ldots, n$, and one can take $g_{k, l}(x, h)$
to be polynomials in $x_1, \ldots, x_n$ and $h_1, \ldots, h_n$ for
each $k$, $l$ as well.

\section{Left translations}
\label{left translations}
\setcounter{equation}{0}

        Let $A$ be an abelian group, and let $N$ be a positive
integer, so that $A^N$ is an abelian group with respect to
coordinatewise addition.  Also let $B(w, z)$ be an $A$-valued function
on $w, z \in A^N$ that is additive in each variable separately, and
let $\diamond$ be the usual group operation on $A^N \times A$
associated to $B(w, z)$.  If $f$ is a function on $A^N \times A$, then put
\begin{equation}
\label{L_{(w, s)}(f)((z, t)) = f((w, s)^{-1} diamond (z, t))}
        L_{(w, s)}(f)((z, t)) = f((w, s)^{-1} \diamond (z, t)).
\end{equation}
Thus $L_{(w, s)}(f)$ is also a function on $A^N \times A$, and
\begin{eqnarray}
\label{L_{(w', s')}(L_{(w, s)}(f))((z, t)) = ...}
        L_{(w', s')}(L_{(w, s)}(f))((z, t))
         & = & (L_{(w, s)}(f))((w', s')^{-1} \diamond (z, t)) \\
 & = & f((w, s)^{-1} \diamond (w', s')^{-1} \diamond (z, t)) \nonumber \\
 & = & f(((w', s') \diamond (w, s))^{-1} \diamond (z, t)) \nonumber \\
 & = & L_{(w', s') \diamond (w, s)}(f)((z, t))   \nonumber
\end{eqnarray}
for every $(w, s), (w', s'), (z, t) \in A^n \times A$.

        Suppose now that $A$ is a commutative ring, and that $B(w, z)$
is given by
\begin{equation}
\label{B(w, z) = sum_j sum_l b_{j, l} w_j z_l, left translations}
        B(w, z) = \sum_{j = 1}^N \sum_{l = 1}^N b_{j, l} \, w_j \, z_l
\end{equation}
for some $b_{j, l} \in A$, $j, l = 1, \ldots, N$.  If $f((z, t))$ is a
formal power series in the commuting indeterminants $z_1, \ldots, z_N$,
$t$ with coefficients in $A$, and if $w_1, \ldots, w_N$, $s$ are
commuting indeterminants that commute with $z_1, \ldots, z_N$, $t$ as
well, then $L_{(w, s)}(f)((z, t))$ may be defined as in (\ref{L_{(w,
s)}(f)((z, t)) = f((w, s)^{-1} diamond (z, t))}), as a formal power
series in all of these variables.  If $w'_1, \ldots, w'_N$, $s'$ are
additional commuting indeterminants that also commute with the
previous indeterminants, then (\ref{L_{(w', s')}(L_{(w, s)}(f))((z,
t)) = ...}) may be interpreted as an identity between formal power
series in all of these variables.

        Similarly, if $f((z, t))$ is a formal polynomial in the
commuting indeterminants $z_1, \ldots, z_N$, $t$, then $L_{(w,
s)}(f)((z, t))$ may be defined as in (\ref{L_{(w, s)}(f)((z, t)) =
f((w, s)^{-1} diamond (z, t))}) as a formal polynomial in $z_1,
\ldots, z_N$, $t$, $w_1, \ldots, w_N$, $s$, and (\ref{L_{(w', s')}(L_{(w,
s)}(f))((z, t)) = ...}) may be interpreted as an identity between
formal polynomials in the appropriate variables.  In this case, we can
also take $(w, s)$ to be an actual element of $A^N \times A$, and
interpret $L_{(w, s)}(f)((z, t))$ as a formal polynomial in $z_1,
\ldots z_N$, $t$, whose coefficients depend on $w_1, \ldots, w_N$, $s$.
If $(w', s')$ is another element of $A^N \times A$, then (\ref{L_{(w',
s')}(L_{(w, s)}(f))((z, t)) = ...}) holds as an identity between
formal polynomials in $z_1, \ldots, z_N$, $t$.  We can also consider the
$A$-valued function on $A^N \times A$ corresponding to $f((z, t))$ as
a formal polynomial in $z_1, \ldots, z_N$, $t$, and then we are back to
the initial formulation of (\ref{L_{(w, s)}(f)((z, t)) = f((w, s)^{-1}
diamond (z, t))}) and (\ref{L_{(w', s')}(L_{(w, s)}(f))((z, t)) = ...}).

\section{$R$-Algebras}
\label{R-algebras}
\setcounter{equation}{0}

        Let $R$ be a commutative ring, and let $A$ be an abelian
group.  Remember that $A$ is said to be an $R$-module if to each $r
\in R$ and $a \in A$ an element $r \, a$ of $A$ such that $a \mapsto r
\, a$ is a homomorphism on $A$ for each $r \in R$, and $r \, (r' \, a)
= (r \, r') \, a$ for every $r, r' \in R$ and $a \in A$.  If $R$ has a
nonzero multiplicative identity element $e$, then it is customary to
ask also that $e \, a = a$ for every $a \in A$.  If $A$ is also
equipped with a binary operation of multiplication so that $A$ is a
ring, then $A$ is said to be an associative $R$-algebra if
multiplication in $A$ is also linear over $R$, in the sense that
\begin{equation}
\label{r (a b) = (r a) b = a (r b)}
        r \, (a \, b) = (r \, a) \, b = a \, (r \, b)
\end{equation}
for every $r \in R$ and $a, b \in A$.  In particular, if $R$ is a
field, then an $R$-module is the same as a vector space over $R$, and
the definition of an $R$-algebra reduces to the more familiar
definition of an algebra over a field.

        If $A$ is a commutative ring that contains $R$ as a subring,
for instance, then $A$ may be considered as an $R$-algebra.  This also
works when $A$ is a nonncommutative ring that contains $R$ as a
central subring, meaning that each element of $R$ commutes with every
element of $A$.  If $A$ is an ideal in $R$, or a quotient of $R$, then
$A$ may be considered as an $R$-algebra in the obvious way.  If $n$ is
a positive integer and $A = R^n$ is equipped with coordinatewise
multiplication, then $A$ is an $R$-algebra with respect to
coordinatewise multiplication by elements of $R$.  Any ring $A$ may be
considered as an algebra over the ring of integers, with $n \, a$
equal to the sum of $n$ $a$'s for each $a \in A$ and $n \in {\bf Z}$.

        If $A$ is an $R$-algebra and $f(x)$ is a formal polynomial in
one variable with coefficients in $R$, then $f$ defines an $A$-valued
function on $A$, by substituting elements of $A$ for $x$.  Similarly,
if $A$ is a commutative $R$-algebra and $f(x)$ is a formal polynomial
in $n$ commuting indeterminants $x_1, \ldots, x_n$ with coefficients
in $R$, then $f$ defines an $A$-valued function on $A^n$, by
substituting elements of $A$ for $x_1, \ldots, x_n$.  As in Section
\ref{polynomials}, sums, products, and compositions of polynomials
correspond to sums, products, and compositions of the corresponding
functions on $A$ or $A^n$, as appropriate.

        Let $N$ be a positive integer, and let $b = \{b_{j, l}\}_{j, l
= 1}^N$ be an $N \times N$ matrix with entries in $R$.  This
determines an $R$-valued function $B(w, z)$ defined for $w, z \in R^N$
in the usual way, with
\begin{equation}
\label{B(w, z) = sum_{j = 1}^N sum_{l = 1}^N b_{j, l} w_j z_k, R-algebras}
        B(w, z) = \sum_{j = 1}^N \sum_{l = 1}^N b_{j, l} \, w_j \, z_k.
\end{equation}
More precisely, $B(w, z)$ may be considered as a formal polynomial in
the commuting indeterminants $w_1, \ldots, w_N$, $z_1, \ldots, z_N$ with
coefficients in $R$, which then defines an $R$-valued function on $R^N
\times R^N$ in the usual way.  If $A$ is a commutative $R$-algebra,
then we also get an $A$-valued function on $A^N \times A^N$, as in the
preceding paragraph.  This can be used to define a group structure on
$A^N \times A$, as before.

\section{Invariant differentiation}
\label{invariant differentiation}
\setcounter{equation}{0}

        Let $R$ be a commutative ring, let $N$ be a positive integer,
and let $w_1, \ldots, w_N$, $s$, $z_1, \ldots, z_N$, $t$, and $h_1,
\ldots, h_N$ be commuting indeterminants.  Also let $b = \{b_{j,
l}\}_{j, l = 1}^N$ be an $N \times N$ matrix with entries in $R$, and
put
\begin{equation}
\label{B(w, z) = sum_{j = 1}^N sum_{l = 1}^N b_{j, l} w_j z_l, differentiation}
        B(w, z) = \sum_{j = 1}^N \sum_{l = 1}^N b_{j, l} \, w_j \, z_l.
\end{equation}
As in the previous section, we may consider $B(w, z)$ as a formal
polynomial in $w_1, \ldots, w_N$, $z_1, \ldots, z_N$ with coefficients
in $R$, which defines an $R$-valued function in the usual way, as well
as a group operation $\diamond$ on $R^N \times R$.

        Let $f((z, t))$ be a formal power series in $z_1, \ldots, z_N$,
$t$ with coefficients in $R$, and consider
\begin{equation}
\label{f((z, t) diamond (h, 0)) = f((z + h, t + B(z, h)))}
        f((z, t) \diamond (h, 0)) = f((z + h, t + B(z, h)))
\end{equation}
as a formal power series in $z_1, \ldots, z_N$, $t$, $h_1, \ldots, h_N$.
If $l = 1, \ldots, N$, then put
\begin{equation}
\label{(D_l(f))((z, t)) = ...}
        (D_l(f))((z, t)) = \frac{\partial f}{\partial z_l}((z, t))
                              + \Big(\sum_{j = 1}^N b_{j, l} \, z_j\Big) \,
                                       \frac{\partial f}{\partial t}((z, t)).
\end{equation}
It is easy to see that
\begin{eqnarray}
\label{f((z, t) diamond (h, 0)) = f((z, t)) + ...}
\lefteqn{f((z, t) \diamond (h, 0))} \\
          & = & f((z, t)) + \sum_{l = 1}^N (D_l(f))((z, t)) \, h_l
                              + \hbox{ higher-order terms}, \nonumber
\end{eqnarray}
where the higher-order terms may be expressed as a sum over $j, l = 1,
\ldots, N$ of $h_j \, h_l$ times formal power series in $z_1, \ldots,
z_N$, $t$, $h_1, \ldots, h_N$.  If $f((z, t))$ is a formal polynomial
in $z_1, \ldots, z_N$, $t$, then (\ref{f((z, t) diamond (h, 0)) = f((z
+ h, t + B(z, h)))}) is a formal polynomial in $z_1, \ldots, z_N$, $t$,
$h_1, \ldots, h_N$, (\ref{(D_l(f))((z, t)) = ...}) is a formal
polynomial in $z_1, \ldots, z_N$, $t$, and the higher-order terms in
(\ref{f((z, t) diamond (h, 0)) = f((z, t)) + ...}) may be expressed as
a sum over $j, l = 1, \ldots, N$ of $h_j \, h_l$ times formal
polynomials in $z_1, \ldots, z_N$, $t$, $h_1, \ldots, h_N$.

        If $f((z, t))$ is a formal power series in $z_1, \ldots, z_N$,
$t$, then
\begin{equation}
\label{L_{(w, s)}(f)((z, t) diamond (h, 0)) = ...}
        L_{(w, s)}(f)((z, t) \diamond (h, 0))
                            = f((w, s)^{-1} \diamond (z, t) \diamond (h, 0))
\end{equation}
may be considered as a formal power series in $w_1, \ldots, w_N$, $s$,
$z_1, \ldots, z_N$, $t$, $h_1, \ldots, h_N$.  If we apply (\ref{f((z,
t) diamond (h, 0)) = f((z, t)) + ...}) to $L_{(w, s)}(f)$, then we get that
\begin{eqnarray}
\label{L_{(w, s)}(f)((z, t) diamond (h, 0)) = L_{w, s}(f)((z, t)) + ...}
\lefteqn{\qquad L_{(w, s)}(f)((z, t) \diamond (h, 0))} \\
 & = & L_{(w, s)}(f)((z, t))
                        + \sum_{l = 1}^N (D_l(L_{(w, s)}(f)))((z, t)) \, h_l
                               + \hbox{ higher-order terms}, \nonumber
\end{eqnarray}
where the higher-order terms may be expressed as a sum over $j, l = 1,
\ldots, N$ of $h_j \, h_l$ times power series in $w_1, \ldots, w_N$,
$s$, $z_1, \ldots, z_N$, $t$, $h_1, \ldots, h_N$.  We can also replace
$(z, t)$ in (\ref{f((z, t) diamond (h, 0)) = f((z, t)) + ...}) with
$(w, s)^{-1} \diamond (z, t)$ to get that
\begin{eqnarray}
\label{f((w, s)^{-1} diamond (z, t) diamond (h, 0)) = ...}
\lefteqn{f((w, s)^{-1} \diamond (z, t) \diamond (h, 0))} \\
 & = & f((w, s)^{-1} \diamond (z, t))
    + \sum_{l = 1}^N (D_l(f))((w, s)^{-1} \diamond (z, t)) \, h_l \nonumber \\
  & &  \qquad + \hbox{ higher-order terms}, \nonumber
\end{eqnarray}
where the higher-order terms may again be expressed as the sum over
$j, l = 1, \ldots, N$ of $h_j \, h_l$ times power series in $w_1,
\ldots, w_N$, $s$, $z_1, \ldots, z_N$, $t$, $h_1, \ldots, h_N$.
Equivalently,
\begin{eqnarray}
\label{L_{(w, s)}(f)((z, t) diamond (h, 0)) = L_{(w, s)}(f)((z, t)) + ..., 2}
\lefteqn{\qquad L_{(w, s)}(f)((z, t) \diamond (h, 0))} \\
 & = & L_{(w, s)}(f)((z, t)) + \sum_{l = 1}^N L_{(w, s)}(D_l(f))((z, t)) \, h_l
                                    + \hbox{ higher-order terms}, \nonumber
\end{eqnarray}
where the higher-order terms are the same as in (\ref{f((w, s)^{-1}
diamond (z, t) diamond (h, 0)) = ...}).

        It follows that
\begin{equation}
\label{D_l(L_{w, s}(f))((z, t)) = L_{(w, s)}(D_l(f))((z, t))}
        D_l(L_{w, s}(f))((z, t)) = L_{(w, s)}(D_l(f))((z, t))
\end{equation}
as power series in $w_1, \ldots, w_N$, $s$, $z_1, \ldots, z_N$, $t$
for $l = 1, \ldots, N$.  If $f((z, t))$ is a polynomial in $z_1,
\ldots, z_N$, $t$, then all of the power series in this discussion may
be taken to be polynomials.  In this case, we can also substitute
specific elements of $R$ for $w_1, \ldots, w_N$, $s$, to get
polynomials in $z_1, \ldots, z_N$, $t$ in (\ref{D_l(L_{w, s}(f))((z,
t)) = L_{(w, s)}(D_l(f))((z, t))}).  When $D_l$ is applied to $L_{(w,
s)}(f)$, note that $w_1, \ldots, w_N$, $s$ are basically constants
with respect to the formal derivatives in $z_1, \ldots, z_N$, $t$.
Similarly, when $L_{(w, s)}$ is applied to $D_l(f)$, note that $L_{(w,
s)}$ is applied to the coefficient of the $\partial f / \partial t$
term, as well as to the derivatives of $f$ with respect to $z_l$ and
$t$ as functions of $z_1, \ldots, z_N$, $t$.

\section{Differentiation and dilations}
\label{differentiation, dilations}
\setcounter{equation}{0}

        Let us continue with the notation and hypotheses in the
previous section.  If $r \in R$, then we have seen that
\begin{equation}
\label{delta_r((z, t)) = (r z, r^2 t), differentiation}
        \delta_r((z, t)) = (r \, z, r^2 \, t)
\end{equation}
defines a homomorphism from $R^N \times R$ into itself with respect to
the group structure corresponding to $\diamond$.  If $f((z, t))$ is a formal
power series in the commuting indeterminants $z_1, \ldots, z_N$, $t$, then
\begin{equation}
\label{delta_r^*(f)((z, t)) = f(delta_r((z, t))) = f((r z, r^2 t))}
        \delta_r^*(f)((z, t)) = f(\delta_r((z, t))) = f((r \, z, r^2 \, t))
\end{equation}
is also a formal power series in $z_1, \ldots, z_N$, $t$.  If $f((z,
t))$ is a polynomial in $z_1, \ldots, z_N$, $t$, then
(\ref{delta_r^*(f)((z, t)) = f(delta_r((z, t))) = f((r z, r^2 t))}) is
also a polynomial in $z_1, \ldots, z_N$, $t$ for each $r \in R$.
Alternatively, one can take $r$ to be another indeterminant that
commutes with $z_1, \ldots, z_N$, $t$, and interpret
(\ref{delta_r^*(f)((z, t)) = f(delta_r((z, t))) = f((r z, r^2 t))}) as
a polynomial or power series in $z_1, \ldots, z_N$, $t$, $r$.  If
$f((z, t))$ is a polynomial, then $f((z, t))$ determines an $R$-valued
function on $R^N \times R$ in the usual way.  In this case,
$\delta_r^*(f)((z, t))$ also determines an $R$-valued function on $R^N
\times R$ for each $r \in R$, which is the same as the function on
$R^N \times R$ associated to $f((z, t))$ composed with the mapping
from $R^N \times R$ into itself associated to $\delta_r$.

        Observe that
\begin{equation}
\label{frac{partial}{partial z_l}(delta_r^*(f))(z, t) = ...}
        \frac{\partial}{\partial z_l}(\delta_r^*(f))(z, t)
         = r \, \frac{\partial f}{\partial z_l}(r \, z, r^2 \, t)
\end{equation}
for $l = 1, \ldots, N$, and
\begin{equation}
\label{frac{partial}{partial t}(delta_r^*(f))(z, t) = ...}
        \frac{\partial}{\partial t}(\delta_r^*(f))(z, t)
         = r^2 \, \frac{\partial f}{\partial t}(r \, z, r^2 \, t)
\end{equation}
for each formal power series $f((z, t))$ in $z_1, \ldots, z_N$, $t$.
Here $r$ may be taken to be any element of $R$, or it may be an
indeterminant that commutes with $z_1, \ldots, z_N$, $t$, so that
(\ref{frac{partial}{partial z_l}(delta_r^*(f))(z, t) = ...}) and
(\ref{frac{partial}{partial t}(delta_r^*(f))(z, t) = ...}) are
interpreted as identities between power series in $z_1, \ldots, z_N$, $t$, $r$.

        If $l = 1, \ldots, N$, then we get that
\begin{equation}
\label{(D_l(delta_r^*(f)))((z, t)) = ...}
 \qquad (D_l(\delta_r^*(f)))((z, t))
            = r \, \frac{\partial f}{\partial z_l}((r \, z, r^2 \, t))
                  + \Big(\sum_{j = 1}^N b_{j, l} \, z_j\Big) \, 
                      r^2 \, \frac{\partial f}{\partial t}((r \, z, r^2 \, t)),
\end{equation}
where $D_l$ is as in (\ref{(D_l(f))((z, t)) = ...}).  Equivalently,
\begin{equation}
\label{(D_l(delta_r^*(f)))((z, t)) = ..., 2}
 \qquad (D_l(\delta_r^*(f)))((z, t))
            = r \, \frac{\partial f}{\partial z_l}((r \, z, r^2 \, t))
                  + r \, \Big(\sum_{j = 1}^N b_{j, l} \, r \, z_j\Big) \,
                        \frac{\partial f}{\partial t}((r \, z, r^2 \, t)),
\end{equation}
which implies that
\begin{equation}
\label{(D_l(delta_r^*(f)))((z, t)) = r delta_r^*(D_l(f))((z, t)}
        (D_l(\delta_r^*(f)))((z, t)) = r \, \delta_r^*(D_l(f))((z, t)).
\end{equation}
Note that $\delta_r^*$ acts on the coefficient of $\partial f /
\partial t$ in $D_l(f)$ on the right side, as well as the derivatives
of $f$ with respect to $z_l$ and $t$ as functions of $z_1, \ldots
z_N$, $t$.  One can also see (\ref{(D_l(delta_r^*(f)))((z, t)) = r
delta_r^*(D_l(f))((z, t)}) by applying (\ref{f((z, t) diamond (h, 0))
= f((z, t)) + ...}) with $z$, $t$, and $h$ replaced with $r \, z$,
$r^2 \, t$, and $r \, h$, and comparing that with (\ref{f((z, t)
diamond (h, 0)) = f((z, t)) + ...})  applied to $\delta_r^*(f)$.
As before, one can take $r$ to be an element of $R$ here, or it may be
an indeterminant that commutes with $z_1, \ldots, z_N$, $t$, so that
(\ref{(D_l(delta_r^*(f)))((z, t)) = r delta_r^*(D_l(f))((z, t)})
is an identity between power series in $z_1, \ldots, z_N$, $t$, $r$.

\section{Real numbers}
\label{real numbers}
\setcounter{equation}{0}

        Let $N$ be a positive integer, and let
\begin{equation}
        B(w, z) = \sum_{j = 1}^N \sum_{l = 1}^N b_{j, l} \, w_j \, z_l
\end{equation}
be a bilinear form on ${\bf R}^N$, so that
\begin{equation}
        (w, s) \diamond (z, t) = (w + z, s + t + B(w, z))
\end{equation}
defines a group structure on ${\bf R}^N \times {\bf R}$, as usual.  Of
course, the real line ${\bf R}$ is a locally compact abelian
topological group with respect to ordinary addition, and Lebesgue
measure on ${\bf R}$ is invariant under translations.  Similarly,
$R^n$ is a locally compact abelian topological group with respect to
coordinatewise addition, and $n$-dimensional Lebesgue measure on ${\bf
R}^n$ is invariant under translations.  It is easy to see that $(N +
1)$-dimensional Lebesgue measure on $R^N \times {\bf R}$ is also
invariant with respect to translations on the left and on the right
corresponding to $\diamond$.

        If $r$ is a real number and $E$ is a Borel subset of ${\bf
R}$, then $r \, E$ is a Borel set in ${\bf R}$ whose Lebesgue measure
is equal to $|r|$ times the Lebesgue measure of $E$, where $|r|$ is
the ordinary absolute value of $r$.  If $E$ is a Borel set in ${\bf
R}^n$, then $r \, E$ is a Borel set in ${\bf R}^n$ with
$n$-dimensional Lebesgue measure equal to $|r|^n$ times the
$n$-dimensional Lebesgue measure of $E$.  As in Section
\ref{dilations} let $\delta_r$ be the parabolic dilation on ${\bf R}^N
\times {\bf R}$ associated to $r \in {\bf R}$, given by
\begin{equation}
        \delta_r((z, t)) = (r \, z, r^2 \, t)
\end{equation}
for each $r, t \in {\bf R}$ and $z \in {\bf R}^N$.  If $E$ is a Borel
set in ${\bf R}^N \times {\bf R}$, then $\delta_r(E)$ is a Borel set
in ${\bf R}^N \times {\bf R}$ for each $r \in {\bf R}$, and the $(N +
1)$-dimensional Lebesgue measure of $\delta_r(E)$ is equal to $|r|^{N
+ 2}$ times the $(N + 1)$-dimensional Lebesgue measure of $E$.

        In this case, there are many natural classes of real and
complex-valued functions on ${\bf R}^N \times {\bf R}$ to consider.
In particular, if $f((z, t))$ is a $C^\infty$ function on ${\bf R}^N
\times {\bf R}$ in the classical sense, then the left translations
$L_{(w, s)}(f)$ of $f$ as in Section \ref{left translations} are
smooth as well, as is the composition of $f$ with a dilation
$\delta_r$ as in the previous paragraph.  In fact, if $f((z, t))$ is a
smooth function on ${\bf R}^N \times {\bf R}$, then $L_{(w, s)}(f)((z,
t))$ is smooth in all of the variables $w_1, \ldots, w_N$, $s$, $z_1,
\ldots, z_N$, $t$, and similarly $f(\delta_r(z, t))$ is smooth in
$z_1, \ldots, z_N$, $t$, $r$.  Note that the space of $C^\infty$
functions on ${\bf R}^N \times {\bf R}$ is a commutative algebra with
respect to pointwise addition and multiplication, which includes
polynomials on ${\bf R}^N \times {\bf R}$ as a subalgebra.

        The formal differential operators discussed in Sections
\ref{invariant differentiation} and \ref{differentiation, dilations}
may be considered as classical differential operators acting on smooth
real or complex-valued functions on ${\bf R}^N \times {\bf R}$.  It is
easy to see that the differential operators in (\ref{(D_l(f))((z, t))
= ...}) have the same invariance properties with respect to
translations and dilations as before.  Of course, formal derivatives
of polynomials or formal power series have many of the same properties
as classical derivatives of smooth functions, such as linearity, the
product rule, and the chain rule.  One can also use Taylor polynomials
or simply the definition of differentiability to analyze the local
behavior of a function at a point, in much the same way that the local
behavior of a polynomial or power series can be analyzed
algebraically.

\section{Another variant}
\label{another variant}
\setcounter{equation}{0}

        Let $A$ and $A'$ be abelian groups, in which the group
operations are expressed by addition, and let $N$ be a positive
integer.  Thus $A^N$ is an abelian group with respect to
coordinatewise addition, as usual.  Also let $B(w, z)$ be an
$A'$-valued function of $w, z \in A^N$ which is additive in each
variable separately, in the sense that
\begin{equation}
        B(w + \widetilde{w}, z) = B(w, z) + B(\widetilde{w}, z)
\end{equation}
and
\begin{equation}
        B(w, z + \widetilde{z}) = B(w, z) + B(w, \widetilde{z})
\end{equation}
for every $w, \widetilde{w}, z, \widetilde{z} \in A^N$.  Under these
conditions,
\begin{equation}
\label{(w, s) diamond (z, t) = (w + z, s + t + B(w, z)), another variant}
        (w, s) \diamond (z, t) = (w + z, s + t + B(w, z))
\end{equation}
defines a group operation on $A^N \times A'$, for the same reasons as before.

        Of course, we have previously discussed this in the case where
$A = A'$.  Even if one starts with $B(w, z)$ as an $A$-valued function
of $w, z \in A^N$ which is additive in each variable, there are
several ways in which one may arrive at $A^N \times A'$ instead of
$A^N \times A$ for some abelian group $A'$.  If $B(w, z)$ happens to
take values in a subgroup $A'$ of $A$, then $A^N \times A'$ is a
subgroup of $A^N \times A$.  Alternatively, if $\phi$ is a
homomorphism from $A$ into $A'$, then
\begin{equation}
        \Phi((z, t)) = (z, \phi(t))
\end{equation}
is a homomorphism from $A^N \times A$ into $A^N \times A'$, using
\begin{equation}
        B'(w, z) = \phi(B(w, z))
\end{equation}
as an $A'$-valued function of $w, z \in A^N$ to define the group
structure on $A^N \times A'$ as before.  In particular, $A$ may simply
be a subgroup of $A'$, so that $A^N \times A$ is a subgroup of $A^N
\times A'$ with respect to the group structure associated to $B(w,
z)$.

        In the situation described at the beginning of the section, it
is easy to see that $\{0\} \times A'$ is a subgroup of $A^N \times A'$
with respect to (\ref{(w, s) diamond (z, t) = (w + z, s + t + B(w,
z)), another variant}), as in Section \ref{another subgroup}.
Moreover,
\begin{equation}
        (w, s) \diamond (0, t) = (0, t) \diamond (w, s) = (w, s + t)
\end{equation}
for every $w \in A^N$ and $s, t \in A'$, and
\begin{equation}
        t \mapsto (0, t)
\end{equation}
defines an isomorphism from $A'$ onto $\{0\} \times A'$.  The natural
projection $\Psi$ from $A^N \times A'$ onto $A^N$ defined by
\begin{equation}
        \Psi((z, t)) = z
\end{equation}
is a homomorphism from $A^N \times A'$ equipped with (\ref{(w, s)
diamond (z, t) = (w + z, s + t + B(w, z)), another variant}) onto
$A^N$ as an abelian group with coordinatewise addition, with kernel
equal to $\{0\} \times A'$.

        Suppose that $C(w, z)$ is another $A'$-valued function of
$w, z \in A^N$ which is additive in each variable separately, and put
\begin{equation}
        \widetilde{B}(w, z) = B(w, z) + C(w, z) + C(z, w).
\end{equation}
This leads to another group operation on $A^N \times A'$, defined by
\begin{equation}
\label{(w, s) widetilde{diamond} (z, t) = ..., another variant}
 (w, s) \widetilde{\diamond} (z, t) = (w + z, s + t + \widetilde{B}(w, z)).
\end{equation}
As in Section \ref{changing variables},
\begin{equation}
\label{(z, t) mapsto (z, t + C(z, z)), another variant}
        (z, t) \mapsto (z, t + C(z, z))
\end{equation}
defines an isomorphism from $A^N \times A'$ with respect to (\ref{(w,
s) diamond (z, t) = (w + z, s + t + B(w, z)), another variant}) onto
$A^N \times A'$ with respect to (\ref{(w, s) widetilde{diamond} (z, t)
= ..., another variant}).  One can also consider $2$-cocycles and
$2$-coboundaries on $A^N$ with values in $A'$, as in Section
\ref{cocycles}.  If a commutative semigroup $\Sigma$ acts on $A$ and
$A'$ in a nice way, then one can let $\Sigma$ act on $A^N \times A'$
as well, as in Section \ref{dilations}.

\section{Some more subgroups}
\label{some more subgroups}
\setcounter{equation}{0}

        Let $A$, $A'$ be abelian groups, let $N$ be a positive
integer, and let $B(w, z)$ be an $A'$-valued function of $w, z \in
A^N$ which is additive in each variable, so that (\ref{(w, s) diamond
(z, t) = (w + z, s + t + B(w, z)), another variant}) defines a group
structure on $A^N \times A'$, as in the previous section.  If $A_1$,
$A_1'$ are subgroups of $A$, $A'$, respectively, such that
\begin{equation}
\label{B(w, z) in A_1'}
        B(w, z) \in A_1'
\end{equation}
for every $w, z \in A_1$, then $A_1^N \times A_1'$ is a subgroup of
$A^N \times A'$ with respect to (\ref{(w, s) diamond (z, t) = (w + z,
s + t + B(w, z)), another variant}), as in Section \ref{subgroups,
quotients}.  If
\begin{equation}
        B(w, z) - B(z, w) \in A_1'
\end{equation}
for every $w \in A^N$ and $z \in A_1^N$, then $A_1^N \times A_1'$ is a
normal subgroup of $A^N \times A'$ with respect to (\ref{(w, s)
diamond (z, t) = (w + z, s + t + B(w, z)), another variant}), as before.
In particular, this holds when
\begin{equation}
        B(w, z), B(z, w) \in A_1'
\end{equation}
for every $w \in A^N$ and $z \in A_1^N$.  In the case where $A_1 = A$,
we get that $A^N \times A_1'$ is a normal subgroup of $A^N \times A'$
with respect to (\ref{(w, s) diamond (z, t) = (w + z, s + t + B(w,
z)), another variant}) as soon as (\ref{B(w, z) in A_1'}) holds for
every $w, z \in A^N$.

        Now let $R$ be a commutative ring, and let $b = \{b_{j,
l}\}_{j, l = 1}^N$ be an $N \times N$ matrix with entries in $R$, so that
\begin{equation}
        B(w, z) = \sum_{j = 1}^N \sum_{l = 1}^N b_{j, l} \, w_j \, z_l
\end{equation}
is an $R$-valued function of $w, z \in R^N$ which is linear over $R$
in each variable separately.  This leads to a group structure
$\diamond$ on $R^N \times R$ in the usual way.  If $R_1$ is a subring
of $R$, then $R_1^N \times R_1$ is a subgroup of $R^N \times R$ with
respect to $\diamond$, as in Section \ref{rings, subrings}.  Let $R_2$
be the subset of $R$ consisting of finite sums of products of at least
two elements of $R_1$.  This is a subring of $R$ which is an ideal in
$R_1$, and $R_2$ is an ideal in $R$ when $R_1$ is an ideal in $R$.  By
construction,
\begin{equation}
\label{B(w, z) in R_2}
        B(w, z) \in R_2
\end{equation}
for every $w, z \in R_1$, and hence $R_1^N \times R_2$ is a subgroup
of $R^N \times R$ with respect to $\diamond$.  More precisely, $R_1^N
\times R_2$ is a normal subgroup of $R_1^N \times R_1$, but it may not
be a normal subgroup of $R^N \times R$, even when $R_1$ is an ideal in
$R$.

        If $r \in R$ and $(z, t) \in R^N \times R$, then put
\begin{equation}
        \delta_r((z, t)) = (r \, z, r^2 \, t),
\end{equation}
as in Section \ref{dilations}.  If $R_1$ is an ideal in $R$, then
\begin{equation}
        \delta_r((z, t)) \in R_1^N \times R_1
\end{equation}
for every $r \in R$ and $(z, t) \in R_1^N \times R_1$, as before.
Similarly,
\begin{equation}
        \delta_r((z, t)) \in R_1^N \times R_2
\end{equation}
for every $r \in R$ and $(z, t) \in R_1^N \times R_2$ in this case.
The same conclusion also holds for every $r \in R_1$ and $(z, t) \in
R^N \times R$ when $R_1$ is an ideal in $R$.

\section{The standard ultranorm}
\label{standard ultranorm}
\setcounter{equation}{0}

        Let $p$ be a prime number, let $n$ be a positive integer, and
consider ${\bf Q}_p^n$ as a vector space over ${\bf Q}_p$ with respect
to coordinatewise addition and scalar multiplication.  The standard
ultranorm on ${\bf Q}_p^n$ is defined by
\begin{equation}
\label{||v||_p = max_{1 le j le n} |v_j|_p}
        \|v\|_p = \max_{1 \le j \le n} |v_j|_p,
\end{equation}
where $v = (v_1, \ldots, v_n) \in {\bf Q}_p^n$.  Thus 
\begin{equation}
        \|r \, v\|_p = |r|_p \, \|v\|_p
\end{equation}
and
\begin{equation}
        \|v + w\|_p \le \max (\|v\|_p, \|w\|_p)
\end{equation}
for every $r \in {\bf Q}_p$ and $v, w \in {\bf Q}_p^n$.  In
particular, this implies that
\begin{equation}
        d_p(v, w) = \|v - w\|_p
\end{equation}
defines an ultrametric on ${\bf Q}_p^n$, which is the standard
ultrametric on ${\bf Q}_p^n$.  Note that the topology on ${\bf Q}_p^n$
determined by the standard ultrametric is the same as the product
topology associated to the usual topology on ${\bf Q}_p$.

        Let $N$ be a positive integer, and take $n = N + 1$ in the
previous discussion, so that ${\bf Q}_p^n$ may be identified with
${\bf Q}_p^N \times {\bf Q}_p$.  Also let
\begin{equation}
        B(w, z) = \sum_{j = 1}^N \sum_{l = 1}^N b_{j, l} \, w_j \, z_l
\end{equation}
be a bilinear form on ${\bf Q}_p^N$, which leads to a group structure
\begin{equation}
        (w, s) \diamond (z, t) = (w + z, s + t + B(w, z))
\end{equation}
on ${\bf Q}_p^N \times {\bf Q}_p$ in the usual way.  Suppose in
addition that $b_{j, l} \in {\bf Z}_p$ for each $j, l = 1, \ldots, N$,
so that ${\bf Z}_p^N \times {\bf Z}_p$ is a subgroup of ${\bf Q}_p^N
\times {\bf Q}_p$ with respect to $\diamond$.  

        As a slight abuse of notation, let us use $\|z\|_p$ for the
standard ultranorm of $z \in {\bf Q}_p^N$, and $\|(z, t)\|_p$ for the
standard ultranorm of $(z, t) \in {\bf Q}_p^N \times {\bf Q}_p \cong
{\bf Q}_p^{N + 1}$.  If $w \in {\bf Z}_p^N$, then it is easy to see
that
\begin{equation}
\label{|B(w, z)|_p le ||z||_p}
        |B(w, z)|_p \le \|z\|_p,
\end{equation}
for every $z \in {\bf Q}_p^N$, and similarly
\begin{equation}
\label{|B(w, z)|_p le ||w||_p}
        |B(w, z)|_p \le \|w\|_p
\end{equation}
for every $w \in {\bf Q}_p^N$ when $z \in {\bf Z}_p^N$.

        It follows that
\begin{equation}
\label{||(w, s) diamond (z, t)||_p le max (||(w, s)||_p, ||(z, t)||_p)}
        \|(w, s) \diamond (z, t)\|_p \le \max (\|(w, s)\|_p, \|(z, t)\|_p)
\end{equation}
for $(w, s), (z, t) \in {\bf Q}_p^N \times {\bf Q}_p$ such that either
$w$ or $z$ is in ${\bf Z}_p$.  If $B(w, z)$ is antisymmetric, then
\begin{equation}
\label{||(z, t)^{-1}||_p = ||(z, t)||_p}
        \|(z, t)^{-1}\|_p = \|(z, t)\|_p
\end{equation}
for every $(z, t) \in {\bf Q}_p^N \times {\bf Q}_p$, because $(z,
t)^{-1} = (-z, -t)$.  Otherwise,
\begin{equation}
        (z, t)^{-1} = (-z, -t + B(z, z)),
\end{equation}
and one can check that (\ref{||(z, t)^{-1}||_p = ||(z, t)||_p}) holds
when $z \in {\bf Z}_p$, since
\begin{equation}
        |B(z, z)|_p \le \|z\|_p.
\end{equation}

        As in Section \ref{subgroups, quotients},
\begin{equation}
 (w, s) \diamond (z, t) \diamond (w, s)^{-1} = (z, t + B(w, z) - B(z, w))
\end{equation}
for every $(w, s), (z, t) \in {\bf Q}_p^N \times {\bf Q}_p$.
If $w \in {\bf Z}_p$, then
\begin{equation}
        |B(w, z)|_p, |B(z, w)|_p \le \|z\|_p,
\end{equation}
by (\ref{|B(w, z)|_p le ||z||_p}) and (\ref{|B(w, z)|_p le ||w||_p}).
Using these observations, one can check that
\begin{equation}
\label{||(w, s) diamond (z, t) diamond (w, s)^{-1}||_p = ||(z, t)||_p}
        \|(w, s) \diamond (z, t) \diamond (w, s)^{-1}\|_p = \|(z, t)\|_p
\end{equation}
for every $(w, s), (z, t) \in {\bf Q}_p^N \times {\bf Q}_p$ such that
$w \in {\bf Z}_p^N$.

        If $(w, s), (z, t) \in {\bf Z}_p^N \times {\bf Z}_p$, then put
\begin{equation}
\label{d((w, s), (z, t)) = ||(w, s)^{-1} diamond (z, t)||_p}
        d((w, s), (z, t)) = \|(w, s)^{-1} \diamond (z, t)\|_p.
\end{equation}
This satisfies the ultrametric version of the triangle inequality,
because of (\ref{||(w, s) diamond (z, t)||_p le max (||(w, s)||_p,
||(z, t)||_p)}).  This is also symmetric in $(w, s)$ and $(z, t)$,
because of (\ref{||(z, t)^{-1}||_p = ||(z, t)||_p}).  Thus (\ref{d((w,
s), (z, t)) = ||(w, s)^{-1} diamond (z, t)||_p}) defines an
ultrametric on ${\bf Z}_p^N \times {\bf Z}_p$, and it is easy to see
that the topology on ${\bf Z}_p^N \times {\bf Z}_p$ determined by
(\ref{d((w, s), (z, t)) = ||(w, s)^{-1} diamond (z, t)||_p}) is the
same as the product topology corresponding to the standard topology on
${\bf Z}_p$.  Note that (\ref{d((w, s), (z, t)) = ||(w, s)^{-1}
diamond (z, t)||_p}) is also invariant under left translations on
${\bf Z}_p^N \times {\bf Z}_p$, by construction.

        Moreover, (\ref{d((w, s), (z, t)) = ||(w, s)^{-1} diamond (z,
t)||_p}) is invariant under translations on ${\bf Z}_p^N \times {\bf
Z}_p$ in the right.  This follows from the invariance of the standard
norm under conjugation by elements of ${\bf Z}_p^N \times {\bf Z}_p$,
as in (\ref{||(w, s) diamond (z, t) diamond (w, s)^{-1}||_p = ||(z,
t)||_p}).  Of course, this ultrametric does not respect dilations on
${\bf Q}_p^N \times {\bf Q}_p$ as in Section \ref{dilations} like the
one in Section \ref{invariant ultrametric}.

\section{A differential equation}
\label{differential equation}
\setcounter{equation}{0}

        Let $R$ be a commutative ring, let $N$ be a positive integer,
and let $b = \{b_{j, l}\}_{j, l = 1}^N$ be an $N \times N$ matrix with
entries in $R$.  Thus
\begin{equation}
        B(w, z) = \sum_{j = 1}^N \sum_{l = 1}^N b_{j, l} \, w_j \, z_l
\end{equation}
is an $R$-valued function of $w, z \in R^N$ which is linear in both
variables, and which leads to a group structure $\diamond$ on $R^N
\times R$ in the usual way.  Also let
\begin{equation}
        \phi(x) = (\phi_1(x), \ldots, \phi_N(x), \phi_{N + 1}(x))
\end{equation}
be an $(N + 1)$-tuple of formal power series in a single variable $x$
with coefficients in $R$, and consider the differential equation
\begin{equation}
\label{phi_{N + 1}'(x) = sum_j sum_l b_{j, l} phi_j(x) phi_l'(x)}
        \phi_{N + 1}'(x) =
          \sum_{j = 1}^N \sum_{l = 1}^N b_{j, l} \, \phi_j(x) \, \phi_l'(x).
\end{equation}
Of course, it is easy to solve this differential equation, at least if
we can divide elements of $R$ by positive integers.

        If $f((z, t))$ is a formal polynomial with coefficients in $R$
in the commuting indeterminants $z_1, \ldots, z_N$, $t$, then
\begin{equation}
\label{f(phi(x)) = f((phi_1(x), ldots, phi_N(x), phi_{N + 1}(x))}
        f(\phi(x)) = f((\phi_1(x), \ldots, \phi_N(x), \phi_{N + 1}(x))
\end{equation}
is a formal power series in $x$ with coefficients in $R$.  This also
works for a formal power series $f((z, t))$ in $z_1, \ldots, z_N$, $t$
when the $\phi_j(x)$'s have constant term equal to $0$ for each $j$.
Under these conditions, the differential equation (\ref{phi_{N +
1}'(x) = sum_j sum_l b_{j, l} phi_j(x) phi_l'(x)}) says exactly that
\begin{equation}
\label{frac{d}{dx} f(phi(x)) = sum_{l = 1}^N phi_l'(x) (D_l(f))(phi(x))}
 \frac{d}{dx} f(\phi(x)) = \sum_{l = 1}^N \phi_l'(x) \, (D_l(f))(\phi(x)),
\end{equation}
where $D_l(f)$ is as in (\ref{(D_l(f))((z, t)) = ...}), by the chain rule.

        Suppose that $(w, s) \in R^N \times R$, and consider
\begin{equation}
\label{psi(x) = (w, s) diamond phi(x)}
        \psi(x) = (w, s) \diamond \phi(x).
\end{equation}
More precisely,
\begin{equation}
        \psi_j(x) = w_j + \phi_j(x)
\end{equation}
for $j = 1, \ldots, N$, and
\begin{equation}
        \psi_{N + 1}(x) = s + \phi_{N + 1}(x)
               + \sum_{j = 1}^N \sum_{l = 1}^N b_{j, l} \, w_j \, \phi_l(x).
\end{equation}
If $\phi(x)$ satisfies (\ref{phi_{N + 1}'(x) = sum_j sum_l b_{j, l}
phi_j(x) phi_l'(x)}), then
\begin{eqnarray}
 \psi_{N + 1}'(x) & = & \phi_{N + 1}'(x) 
          + \sum_{j = 1}^N \sum_{l = 1}^N b_{j, l} \, w_j \, \phi_l'(x) \\
 & = & \sum_{j = 1}^N \sum_{l = 1}^N b_{j, l} \, \phi_j(x) \, \phi_l'(x)
    + \sum_{j = 1}^N \sum_{l = 1}^N b_{j, l} \, w_j \, \phi_l'(x) \nonumber \\
 & = & \sum_{j = 1}^N \sum_{l = 1}^N b_{j, l} \, (w_j + \phi_j(x))
                                                     \, \phi_l'(x) \nonumber \\
 & = & \sum_{j = 1}^N \sum_{l = 1}^N b_{j, l} \, \psi_j(x) \, \psi_l'(x).
                                                                   \nonumber
\end{eqnarray}
This shows that $\phi(x)$ also satisfies (\ref{phi_{N + 1}'(x) = sum_j
sum_l b_{j, l} phi_j(x) phi_l'(x)}) when $\phi(x)$ does.

        If $R$ is the field ${\bf R}$ of real numbers, then
(\ref{phi_{N + 1}'(x) = sum_j sum_l b_{j, l} phi_j(x) phi_l'(x)})
makes sense for any differentiable function defined on an interval in
the real line, with values in ${\bf R}^N \times {\bf R}$.  Similarly,
(\ref{frac{d}{dx} f(phi(x)) = sum_{l = 1}^N phi_l'(x)
(D_l(f))(phi(x))}) follows from the chain rule when $f((z, t))$ is a
differentiable real-valued function on ${\bf R}^N \times {\bf R}$ and
$\phi$ satisfies (\ref{phi_{N + 1}'(x) = sum_j sum_l b_{j, l} phi_j(x)
phi_l'(x)}).  If $\phi(x)$ is a differentiable path in ${\bf R}^n
\times {\bf R}$ that satisfies (\ref{phi_{N + 1}'(x) = sum_j sum_l
b_{j, l} phi_j(x) phi_l'(x)}), then the computation in the previous
paragraph shows that any path $\psi(x)$ obtained by translating
$\phi(x)$ on the left in ${\bf R}^N \times {\bf R}$ with respect to
$\diamond$ also satisfies (\ref{phi_{N + 1}'(x) = sum_j sum_l b_{j, l}
phi_j(x) phi_l'(x)}).

\section{Power series on ${\bf Q}_p$}
\label{power series on Q_p}
\setcounter{equation}{0}

        Let $p$ be a prime number, let $k$ be an integer, and let
$\{a_j\}_{j = 0}^\infty$ be a sequence of $p$-adic numbers such that
\begin{equation}
        \lim_{j \to \infty} |a_j|_p \, p^{-j k} = 0.
\end{equation}
This implies that
\begin{equation}
        f(x) = \sum_{j = 0}^\infty a_j \, x^j
\end{equation}
converges for each $x$ in the set
\begin{equation}
\label{{x in {bf Q}_p : |x|_p le p^{-k}}}
        \{x \in {\bf Q}_p : |x|_p \le p^{-k}\}.
\end{equation}
More precisely, it is easy to see that the partial sums of this power
series converge uniformly to $f(x)$ on (\ref{{x in {bf Q}_p : |x|_p le
p^{-k}}}).  In particular, this implies that $f(x)$ is continuous as a
${\bf Q}_p$-valued function on (\ref{{x in {bf Q}_p : |x|_p le
p^{-k}}}), since polynomials are continuous.

        Let $\{b_j\}_{j = 1}^\infty$ be another sequence of $p$-adic
numbers such that
\begin{equation}
        \lim_{j \to \infty} |b_j|_p \, p^{-j k} = 0,
\end{equation}
and put
\begin{equation}
        g(x) = \sum_{j = 0}^\infty b_j \, x^j.
\end{equation}
Of course,
\begin{equation}
        f(x) + g(x) = \sum_{j = 0}^\infty (a_j + b_j) \, x^j
\end{equation}
for each $x$ in (\ref{{x in {bf Q}_p : |x|_p le p^{-k}}}).  Similarly,
\begin{equation}
\label{f(x) g(x) = sum_{n = 0}^infty c_n x^n}
        f(x) \, g(x) = \sum_{n = 0}^\infty c_n \, x^n
\end{equation}
for every $x$ in (\ref{{x in {bf Q}_p : |x|_p le p^{-k}}}), where
\begin{equation}
        c_n = \sum_{j = 0}^n a_j \, b_{n - j}.
\end{equation}
Note that
\begin{equation}
        |c_n|_p \, p^{-k n} \le \max_{1 \le j \le n} (|a_j|_p \, p^{-j k}) \,
                                            (|b_{n - j}|_p \, p^{-k (n - j)})
\end{equation}
for each $n$, which implies that
\begin{equation}
        \lim_{n \to \infty} |c_n|_p \, p^{-k n} = 0,
\end{equation}
so that the right side of (\ref{f(x) g(x) = sum_{n = 0}^infty c_n
x^n}) converges on (\ref{{x in {bf Q}_p : |x|_p le p^{-k}}}).  In
order to check that (\ref{f(x) g(x) = sum_{n = 0}^infty c_n x^n})
holds on (\ref{{x in {bf Q}_p : |x|_p le p^{-k}}}), one can
approximate $f(x)$ and $g(x)$ by finite sums.

        The corresponding power series for the derivative,
\begin{equation}
\label{f'(x) = sum_{j = 1}^infty j a_j x^{j - 1}, power series on Q_p}
        f'(x) = \sum_{j = 1}^\infty j \, a_j \, x^{j - 1},
\end{equation}
also converges when the series for $f(x)$ converges, because $|j|_p
\le 1$ for every integer $j$.  Using standard arguments, one can verify that
\begin{equation}
        \lim_{h \to 0} \frac{f(x + h) - f(x)}{h} = f'(x)
\end{equation}
for every $x \in {\bf Q}_p$ such that $|x|_p \le p^{-k}$.  Of course,
the derivative is linear in $f$, and satisfies the usual product rule.

        Suppose now that we take $R = {\bf Q}_p$ in the previous
section.  In this case, we can take $\phi(x)$ to be an $(N + 1)$-tuple
of convergent power series on the set (\ref{{x in {bf Q}_p : |x|_p le
p^{-k}}}), and consider the same differential equation (\ref{phi_{N +
1}'(x) = sum_j sum_l b_{j, l} phi_j(x) phi_l'(x)}) as before.  If
$f((z, t))$ is a polynomial in $z_1, \ldots, z_N$, $t$ with
coefficients in ${\bf Q}_p$, then $f(\phi(x))$ is given by a
convergent power series on (\ref{{x in {bf Q}_p : |x|_p le p^{-k}}}),
and we still have (\ref{frac{d}{dx} f(phi(x)) = sum_{l = 1}^N
phi_l'(x) (D_l(f))(phi(x))}) when $\phi(x)$ satisfies (\ref{phi_{N +
1}'(x) = sum_j sum_l b_{j, l} phi_j(x) phi_l'(x)}).  If $\psi(x)$ is
as in (\ref{psi(x) = (w, s) diamond phi(x)}), then $\psi(x)$ is also
given by convergent power series on (\ref{{x in {bf Q}_p : |x|_p le
p^{-k}}}), and satisfies (\ref{phi_{N + 1}'(x) = sum_j sum_l b_{j, l}
phi_j(x) phi_l'(x)}) when $\phi(x)$ does, for the same reasons as before.

\section{Power series on ${\bf Q}_p^n$}
\label{power series on Q_p^n}
\setcounter{equation}{0}

        Let $p$ be a prime number, let $n$ be a positive integer, and
suppose that for each multi-index $\alpha = (\alpha_1, \ldots,
\alpha_n)$ we have a $p$-adic number $A_\alpha$ such that
\begin{equation}
\label{lim_{|alpha| to infty} |A_alpha|_p = 0}
        \lim_{|\alpha| \to \infty} |A_\alpha|_p = 0.
\end{equation}
In particular, this implies that
\begin{equation}
\label{lim_{L to infty} (sum_{|alpha| = L} A_alpha) = 0}
        \lim_{L \to \infty} \Big(\sum_{|\alpha| = L} A_\alpha \Big) = 0
\end{equation}
in ${\bf Q}_p$, and hence that
\begin{equation}
\label{sum_{L = 0}^infty (sum_{|alpha| = L} A_alpha)}
        \sum_{L = 0}^\infty \Big(\sum_{|\alpha| = L} A_\alpha\Big)
\end{equation}
converges in ${\bf Q}_p$.  We can use this last sum as a definition of
\begin{equation}
\label{sum_alpha A_alpha}
        \sum_\alpha A_\alpha,
\end{equation}
where the sum is taken over all multi-indices $\alpha$, but in fact
this sum converges in a much stronger sense, that does not depend on
the way that the terms are arranged, because of (\ref{lim_{|alpha| to
infty} |A_alpha|_p = 0}) and the ultrametric version of the triangle
inequality.

        Let $k = (k_1, \ldots, k_n)$ be an $n$-tuple of integers, and suppose
that for each multi-index $\alpha$ we have $a_\alpha \in {\bf Q}_p$ such that
\begin{equation}
\label{lim_{|alpha| to infty} |a_alpha|_p p^{- alpha cdot k} = 0}
        \lim_{|\alpha| \to \infty} |a_\alpha|_p \, p^{- \alpha \cdot k} = 0,
\end{equation}
where $\alpha \cdot k = \sum_{j = 1}^n \alpha_j \, k_j$.  This implies that
\begin{equation}
\label{f(x) = sum_alpha a_alpha x^alpha, power series on Q_p^n}
        f(x) = \sum_\alpha a_\alpha \, x^\alpha
\end{equation}
converges in ${\bf Q}_p$ as in the previous paragraph for each $x$ in the set
\begin{equation}
\label{{x in {bf Q}_p^n : |x_j|_p le p^{-alpha_j k_j} for j = 1, ldots, n}}
        \{x \in {\bf Q}_p^n : |x_j|_p \le p^{-\alpha_j \, k_j}
                                  \hbox{ for } j = 1, \ldots, n\}.
\end{equation}
The partial sums of (\ref{f(x) = sum_alpha a_alpha x^alpha, power
series on Q_p^n}) converge uniformly to $f(x)$ on this set, so that
$f(x)$ is a continuous ${\bf Q}_p$-valued function on (\ref{{x in {bf
Q}_p^n : |x_j|_p le p^{-alpha_j k_j} for j = 1, ldots, n}}).  As
usual, sums and products of convergent power series on (\ref{{x in {bf
Q}_p^n : |x_j|_p le p^{-alpha_j k_j} for j = 1, ldots, n}}) correspond
to convergent power series on this set as well.  Similarly, formal
derivatives of a convergent power series on this set also converge on
this set, and represent the derivatives of the power series in the
classical sense of limits of difference quotients.

        Now let $N$ be a positive integer, and take $n = N + 1$ in the
preceding discussion, with ${\bf Q}_p^n$ identified with ${\bf Q}_p^N
\times {\bf Q}_p$.  Let
\begin{equation}
        B(w, z) = \sum_{j = 1}^N \sum_{l = 1}^N b_{j, l} \, w_j \, z_l
\end{equation}
be a bilinear form on ${\bf Q}_p^N$, with coefficients $b_{j, l} \in
{\bf Z}_p$ for each $j$, $l$, and let $\diamond$ be the corresponding
group structure on ${\bf Q}_p^N \times {\bf Q}_p$.  Also let $k$ be an
integer, and consider
\begin{equation}
\label{{(z, t) in Q_p^N times Q_p : |z_j|_p le p^{-k}, |t|_p le p^{-2 k}}}
 \quad  \{(z, t) \in {\bf Q}_p^N \times {\bf Q}_p : |z_j|_p \le p^{-k}
              \hbox{ for } j = 1, \ldots, N, \hbox{ and } |t|_p \le p^{-2 k}\}.
\end{equation}
This is the version of (\ref{{x in {bf Q}_p^n : |x_j|_p le p^{-alpha_j
k_j} for j = 1, ldots, n}}) that we shall use in this situation, and which
is a subgroup of ${\bf Q}_p^N \times {\bf Q}_p$.

        Suppose that $f((z, t))$ is defined by a convergent power
series on (\ref{{(z, t) in Q_p^N times Q_p : |z_j|_p le p^{-k}, |t|_p
le p^{-2 k}}}), let $(w, s)$ be an element of (\ref{{(z, t) in Q_p^N
times Q_p : |z_j|_p le p^{-k}, |t|_p le p^{-2 k}}}), and consider
\begin{equation}
\label{L_{(w, s)}(f)((z, t)) = f((w, s)^{-1} diamond (z, t)), power series}
        L_{(w, s)}(f)((z, t)) = f((w, s)^{-1} \diamond (z, t)),
\end{equation}
as in Section \ref{left translations}.  This is also a function on
(\ref{{(z, t) in Q_p^N times Q_p : |z_j|_p le p^{-k}, |t|_p le p^{-2
k}}}), since (\ref{{(z, t) in Q_p^N times Q_p : |z_j|_p le p^{-k},
|t|_p le p^{-2 k}}}) is a subgroup of ${\bf Q}_p^N \times {\bf Q}_p$
with respect to $\diamond$.  One can check that the power series
corresponding to (\ref{L_{(w, s)}(f)((z, t)) = f((w, s)^{-1} diamond
(z, t)), power series}) converges on (\ref{{(z, t) in Q_p^N times Q_p
: |z_j|_p le p^{-k}, |t|_p le p^{-2 k}}}), and represents (\ref{L_{(w,
s)}(f)((z, t)) = f((w, s)^{-1} diamond (z, t)), power series}) on this set.
In addition, one can take invariant derivatives of $f((z, t))$, with
properties as in Sections \ref{invariant differentiation} and
\ref{differentiation, dilations}.  Note that the usual dilations
$\delta_r((z, t)) = (r \, z, r^2 \, t)$ map (\ref{{(z, t) in Q_p^N
times Q_p : |z_j|_p le p^{-k}, |t|_p le p^{-2 k}}}) into itself
when $r \in {\bf Z}_p$.

\end{document}